\documentclass[11pt]{amsart}
\usepackage{amsmath,amssymb,amsfonts,url,mathptmx}
\numberwithin{equation}{section}
\newtheorem{Def}{Definition}[section]
\newtheorem{thm}{Theorem}[section]
\newtheorem{lem}{Lemma}[section]
\newtheorem{rem}{Remark}[section]

\newcommand{\hdot}{^\text{\r{}}\hspace{-.33cm}H}

\begin{document}
\title[Liouville system]{On Liouville systems at critical parameters, Part 1: one bubble} 

\author{Chang-shou Lin}
\address{Department of Mathematics\\
        Taida Institute of Mathematical Sciences\\
        National Taiwan University\\
         Taipei 106, Taiwan } \email{cslin@math.ntu.edu.tw}

\author{Lei Zhang}
\address{Department of Mathematics\\
        University of Florida\\
        358 Little Hall P.O.Box 118105\\
        Gainesville FL 32611-8105}
\email{leizhang@ufl.edu}
\thanks{Zhang is supported in part by NSF Grant 0900864 (1027628)}

\date{\today}

%%%%%%%%%%%%%%%%%%%%%%%%%%%%%%%%%%%%%%%%%%%%%
\begin{abstract}
In this paper we consider bubbling solutions to the general Liouville system:
\begin{equation}\label{abeq1}
\Delta_g u_i^k+\sum_{j=1}^n a_{ij}\rho_j^k\left(\frac{h_j e^{u_j^k}}{\int h_j e^{u_j^k}}-1\right)=0\quad\text{in \,}M,\,\, i=1,...,n \,\, (n\ge 2)
\end{equation}
where $(M,g)$ is a Riemann surface,  and $A=(a_{ij})_{n\times n}$ is a
constant non-negative matrix and $\rho_j^k\to \rho_j$ as $k\to \infty$. Among other things we prove the following sharp estimates.
\begin{enumerate}
\item The location of the blowup point.
\item The convergence rate of $\rho_j^k-\rho_j$, $j=1,..,n$.
\end{enumerate}
These results are of fundamental importance for constructing bubbling solutions. It is interesting to compare the difference
between the general Liouville system and the $SU(3)$ Toda system on estimates (1) and (2).
\end{abstract}
%%%%%%%%%%%%%%%%%%%%%%%%%%%%%%%%%%%%%%%%%%%%%

\maketitle

\section{Introduction}
Let $(M,g)$ be a compact Riemann surface whose volume is normalized to be $1$, $h_1,...,h_n$ be positive $C^3$ functions on $M$, $\rho_1,..,\rho_n$ be nonnegative
constants. In this article we continue our study of the following Liouville system defined on $(M,g)$:
\begin{equation}\label{mainsys}
\Delta_g u_i+\sum_{j=1}^n\rho_j a_{ij} (\frac{h_je^{u_j}}{\int_Mh_je^{u_j}dV_g}-1)=0,\quad i\in I:=\{1,..,n\}
\end{equation}
where $dV_g$ is the volume form, $A=(a_{ij})$ is a non-negative constant matrix, $\Delta_g$ is the Laplace-Beltrami operator
($-\Delta_g\ge 0$).
When $n=1$ and $a_{11}=1$, equation \eqref{mainsys} is the mean field equation of the Liouville type:
\begin{equation}\label{equfromsys}
\Delta_g u+\rho\left(\frac{h e^u}{\int_M h e^udV_g}-1\right)=0\quad\text{in \,}M.
\end{equation}
Therefore, the Liouville system \eqref{mainsys} is a natural extension of the classical Liouville equation,
which has been extensively studied for the past three decades. Both the Liouville equation and the Liouville system are related to various fields of geometry, Physics, Chemistry and Ecology. For example
in conformal geometry, when $\rho=8\pi$ and $M$ is the sphere $\mathbb S^2$, equation \eqref{equfromsys} is equivalent to the famous Nirenberg problem.
For a bounded domain in $\mathbb R^2$ and $n=1$, a variant of \eqref{equfromsys} can be derived from the mean field limit of Euler flows or spherical Onsager vortex theory,
as studied by Caglioti, Lions, Marchioro and Pulvirenti\,\cite{lion1, lion2}, Kiessling\,\cite{kiess}, Chanillo and Kiessling\,\cite{chanillo1} and Lin \cite{linarch}.
In classical gauge field theory, equation \eqref{mainsys} is closely related to the Chern-Simons-Higgs equation for the abelian case,
see \cite{caffarelli,hong,jackiw,yang}.
Various Liouville systems are also used to describe models in the theory of self-gravitating systems \cite{aly}, Chemotaxis\,\cite{childress, keller}, in the physics of charged particle beams\,\cite{bennet,debye,kiess2,kiessling2},
in the non-abelian Chern-Simons-Higgs theory\,\cite{dunne,jostlinwang,yang} and other gauge field models \cite{dziar,phys,kimleelee}.
For recent developments of these subjects or related Liouville systems in more general
settings, we refer the readers to \cite{barto1, barto3, ChenLin1, chenlin2, chenlinnew, clwang, CSW, CSW1, licmp, linduke, linarch, lwy,lwz, linzhang1, nolasco2, rubinstein, spruck, wolansky1, wolansky2, zhangcmp, zhangccm} and the references therein.

For any solution $u$ of \eqref{equfromsys}, clearly adding any constant to $u$ gives another solution. So it is nature to assume $u\in \,\hdot^1(M)$, where
\begin{equation*}
\hdot^1(M)=\left\{
u\in L^2(M)\,\Big|\, |\nabla_g u|\in L^2(M)\text{ and }\int_M  u\,dV_g=0
\right\}.
\end{equation*}
Corresponding to (\ref{mainsys}) we set
$$\hdot^{1,n}=\,\hdot^1(M)\times\cdots\times\,\hdot^1(M) $$ to be the space for solutions.
For any $\rho=(\rho_1,\cdots, \rho_n)$, $\rho_i>0 (i\in I=\{1,...,n\})$, let $\varPhi_\rho$ be a nonlinear functional defined in \, $\hdot^{1,n}$ by
$$
\varPhi_\rho(u)=\frac{1}{2}\sum_{i,j\in I} a^{ij}\int_M \nabla_g u_i\cdot \nabla_g u_j dV_g-\sum_{j\in I}\rho_j\log\int_M h_j e^{u_j}dV_g
$$
where $(a^{ij})_{n\times n}$ is the inverse of $A=(a_{ij})_{n\times n}$.
It is easy to see that equation \eqref{mainsys} is the Euler-Lagrangian equation of $\varPhi_\rho$.

If the matrix $A$ satisfies the following two assumptions:
\begin{eqnarray*}
&& (H1):\quad \mbox{ $A$ is symmetric, nonnegative, irreducible and invertible}. \\
&&(H2): \quad a^{ii}\leq 0,\,\, \forall i\in I, \quad a^{ij}\geq 0, \,\, \forall i\ne j,
\quad \sum_{j\in I}a^{ij}\geq 0, \,\, \forall i\in I,
\end{eqnarray*}
the authors prove in \cite{linzhang2} that for $\rho$ satisfying
\begin{equation}\label{gammaN}
8\pi N\sum_{i\in I}\rho_i<\sum_{i,j\in I}a_{ij}\rho_i\rho_j<8\pi (N+1)\sum_{i\in I}\rho_i,
\end{equation}
there is a priori estimate for all solutions $u$ to (\ref{mainsys}),
and the Leray-Schauder degree $d_{\rho}$ for equation (\ref{mainsys}) is
$$
d_{\rho}=
\frac{1}{N!}\bigg ((-\chi_M+1)...(-\chi_M+N) \bigg )\quad \text{if } \mbox{ (\ref{gammaN}) holds }
$$
where $\chi_M$ is the Euler characteristic of $M$. Moreover, if $\rho^k$ tends to the hyper-surface
$\{\rho;\quad 8\pi N\sum_{i\in I}\rho_i=\sum_{i,j\in I}a_{ij}\rho_i\rho_j\}$, there exist exactly $N$ disjoint blowup points (see \cite{linzhang2}).

The proof of the a priori bound in \cite{linzhang2} relies on the sharp estimate for a sequence of bubbling solutions to (\ref{mainsys}).
Let $u^k$ be the blowup solutions corresponding to $\rho^k$ and  $B(p_t,\delta_0)$ ($t=1,..,N$) be disjoint balls around distinct blowup points in $M$. Then under assumptions $(H1)$ and $(H2)$,
the behavior of $u^k$ around any $p_t$ is fully bubbling, that is, the maximum values of any components of $u^k$ in any of the balls are of the same magnitude:
$$\max_{B(p_t,\delta_0)} u_i^k = \max_{B(p_t,\delta_0)} u_j^k +O(1), \quad \forall i, j\in I. $$
Moreover, after a suitable scaling around each blowup point $p_t$, $u^k$ converges to an entire solution
$\tilde U=(\tilde U_1,..,\tilde U_n)$ of the following Liouville system:
\begin{equation}\label{12apr23e2}
\left\{\begin{array}{ll}
\Delta \tilde U_i+\sum_{j=1}^n a_{ij}e^{\tilde U_j}=0,\quad \mbox{in }\,\, \mathbb R^2, \\
\\
\int_{\mathbb R^2}e^{\tilde U_i}<\infty, \quad \tilde U_i \mbox{ is radial},\,\, \forall i\in I.
\end{array}
\right.
\end{equation}

One may expect the limiting entire solution to be different around each blowup point, however the authors proved that
$\tilde U$ is independent of blowup points, and only depends on the ratio of $\rho_1^k-\rho_1:\rho_2^k-\rho_2:..\rho_n^k-\rho_n$ (see \cite{linzhang1}).  Naturally it leads to the question: how to construct bubbling solutions with the help of this information?

In this paper and subsequent ones, we are devoted to study the bubbling phenomenon of Liouville systems: how to accurately estimate the bubbling solutions of (\ref{mainsys}) and how to construct them. These are quite challenging analytic problems. In general, blowup analysis for a system of equations is much harder than that for the single equation. One reason is that the Pohozaev identity, a balancing condition, is no longer so powerful as in the scalar case. Another reason is that there are too many entire solutions: the parameter
$\sigma=(\sigma_1,...,\sigma_n)$ ($\displaystyle{\sigma_i=\frac 1{2\pi}\int_{\mathbb R^2}e^{u_i}}$), which represents the energy of the entire solution,
 forms a submanifold of $n-1$ dimension. However, for the Liouville equation, the energy is just one number: $\int_{\mathbb R^2}e^u=8\pi$.

 In this article we consider the case of one blowup point, and always assume $(H1)$ only. Let $\rho=(\rho_1,...,\rho_n)$ and
 \begin{equation}\label{12jun20e1}
 \Lambda_J(\rho)=8\pi \sum_{i\in J} \rho_i-\sum_{i,j\in J}a_{ij}\rho_i\rho_j
 \end{equation}
 for any $J\subset I:=\{1,...,n\}$. Define
 $$\Gamma_1=\{\rho; \,\, \Lambda_I(\rho)=0\quad \mbox{ and }\Lambda_J(\rho)>0 \mbox{ for all } \emptyset \subsetneq J \subsetneq I \,\, \}. $$
 Note that if $(H2)$ also holds also, then $\Lambda_I(\rho)=0$ implies $\Lambda_J(\rho)>0$ for all $J\subsetneq I$ ( see \cite{linzhang2}). For any $\rho$ we define $(m_1,...,m_n)$ by
 \begin{equation}\label{12apr23e1}
 m_i=\frac 1{2\pi} \sum_{j\in I} a_{ij}\rho_j.
 \end{equation}
 The quantity $m_i$ can be interpreted by the entire solution $\tilde U$ of (\ref{12apr23e2}). In fact
 $$\rho_i=\int_{\mathbb R^2} e^{\tilde U_i}dx, \quad i\in I, $$
 and
 \begin{equation}\label{12apr23e3}
 \tilde U_i(x)=-m_i \log |x| +O(1),\quad \mbox{ for $|x|$ near infinity. }
 \end{equation}
 The integrability of $e^{\tilde U_i}$ implies $m_i>2$ for all $i$.
 On the other hand $\Lambda_I(\rho)=0$ can be written as $\sum_{i\in I}(m_i-4)\rho_i=0$. Thus either $\min\{m_1,...,m_n\}<4$ or $m_i=4$ for all $i\in I$.
 We also note that (\ref{12apr23e3}) implies that $\Gamma_1$ is a smooth submanifold because the normal vector at $\rho\in \Gamma_1$:
 $$(\sum_{j\in I} a_{ij} \rho_j-4\pi,...,\sum_{j\in I} a_{nj}\rho_j -4\pi), $$
has all its components positive.

 The asymptotic behavior of $\tilde U_i(x)$ shows that the decay rate of $e^{\tilde U_i(x)}$ is $O(|x|^{-m})$, where
 $$m=\min\{m_1,...,m_n\}. $$
 In this article we define $Q\in \Gamma_1$ with $m=4$, i.e. $m_i=4$ for all $i$. Thus the decay rate of $e^{\tilde U_i}$ for $\rho=Q$ is $O(|x|^{-4})$. The difference on the decay rate for $Q$ and $\rho\neq Q$ will have great effects on bubbling analysis later.

 Let $u^k=(u_1^k,...,u_n^k)$ be a sequence of blow up solutions to (\ref{mainsys}) with $\rho=\rho^k$ such that $\rho^k\to \rho\in \Gamma_1$. The point $Q$ defined above is of particular importance, the readers will see that in our main theorems below, the asymptotic behavior of blowup solutions, the nature of $\Lambda_I(\rho^k)$ and the location of blowup point are all significantly different depending on $\rho=Q$ or not.

  Let $p$ be the blowup point of $u^k$ and
 \begin{equation}\label{12apr23e4}
 M_k=\max_{B(p,\delta)} \bigg ( u_1^k(x)-\log \int_M h_1e^{u_1^k}dV_g \bigg ),
 \end{equation}
 \begin{equation}\label{12apr23e5}
 \epsilon_k=e^{-\frac 12 M_k}.
 \end{equation}
 Since $\rho^k\to \Gamma_1$, there is only one blowup point $p$. It is easy to see that $u^k$ fully blows up at $p$ (see Lemma \ref{ful}):
 \begin{equation}\label{12apr23e6}
 \max_{B(p, \delta_0)}(u_i^k(x)-\log \int_M h_i e^{u_i^k}dV_g)=M_k+O(1), \quad i\in I.
 \end{equation}

 Our first result is on the location of the blowup point $p$. Let $p_k\to p$ be where the maximum of $\{u_1^k,..,u_n^k\}$ is attained, then we have
\begin{thm}\label{locations}
Let $\rho^k\to \rho\in \Gamma_1$ and all $\rho_i^k-\rho_i$ have the same sign.
\begin{enumerate}
\item If $\rho\neq Q$, then
\begin{equation}\label{11july13e7}
\sum_{i\in I} \bigg (\nabla (\log h_i)(p_k)+2\pi m\nabla_1 \gamma(p_k,p_k)\bigg )
\rho_{i} =O(\epsilon_k^{m-2}).
\end{equation}
\item
 If $\rho=Q$, then
\begin{equation}\label{11july13e8}
\sum_i\bigg (\nabla (\log h_i)(p_k)+8\pi  \nabla_1 \gamma(p_k,p_k) \bigg )\rho_{i}=O(\epsilon_k^2\log \epsilon_k^{-1}).
\end{equation}
where $\nabla_1$ denotes the derivative with respect to the first variable, and $\gamma(x,y)$ stands for the regular part of the Green's function.
\end{enumerate}
\end{thm}

\medskip

Our second result is about the decay rate of $\Lambda_I(\rho^k)$. To state the result, we need to define the following quantity $D_i$ ($i\in I=\{1,...,n\}$):
\begin{equation}\label{leadingD}
D_i=\lim_{\delta_0\to 0}\big (\delta_0^{2-m}-\frac{m-2}{2\pi}
\int_{M\setminus B(p,\delta_0)}\frac{h_i(x)}{h_i(p)}e^{2\pi m (G(x,p)-\gamma(p,p))}dV_g \big )
\end{equation}
provided that $m<4$. The limit is well defined if $m<4$, see section 6.
\begin{thm}\label{rhoto1mle4} Suppose $\rho^k\to \rho\in \Gamma_1$ and $\rho\neq Q$, if all $\rho_i^k-\rho_i$ have the same sign, then
$$
\Lambda_I(\rho^k)
=8\pi^2\sum_{i\in I_1}(e^{c_i}D_i+o(1))\epsilon_k^{m-2},
$$
where $I_1$ is a subset of $I$ where $m_i=m$ for all $i\in I_1$, $c_i$ is a constant determined in (\ref{12jun29e1}). $o(1)\to 0$
as $k\to \infty$.
\end{thm}
If $M$ is a flat torus with fundamental cell domain
$\Omega\subset \mathbb R^2$, then $D_i$ can be written as
$$D_i=\frac{m-2}{2\pi }\bigg (
\int_{\mathbb R^2\setminus \Omega} \frac{1}{|x-p|^m}dx-\lim_{\delta_0\to 0}\int_{\Omega\setminus B_{\delta_0}}
\frac{H_i(x,p)}{|x-p|^m}dx\bigg ) $$
where
$$H_i(x,p)=\frac{h_i(x)}{h_i(p)}e^{2\pi m(\gamma(x,p)-\gamma(p,p))}-1, \quad i\in I. $$
See \cite{ccl} and \cite{clwang} for related discussions.
\begin{rem} The assumption that all $\rho_i^k-\rho_i$ have the same sign seems due to some technical difficulties. When $n=2$
and $\rho\neq Q$, this assumption is not needed for both Theorem \ref{locations} and Theorem \ref{rhoto1mle4}.
\end{rem}

\begin{thm}\label{rhotome4}
Suppose $\rho^k\to \rho\in \Gamma_1$ and $\rho=Q$. If all $\rho_i^k-Q_i$ have the same sign, then
\begin{equation}\label{12apr23e10}
\Lambda_I(\rho^k)
=-16\pi^2(\sum_{i\in I} b_i e^{c_i}+o(1))\epsilon_k^2\log \epsilon_k^{-1}.
\end{equation}
where
$$
b_i=\frac 14\bigg (\Delta \, \log h_i(p)-2K(p)+8\pi+|\nabla \, \log h_i(p)+8\pi \nabla_1\gamma(p,p)|^2\bigg ),
$$
$c_i$ is determined in (\ref{12jun29e1}).
\end{thm}

Important information on bubbling solutions can be observed on the two cases:
$\rho\neq Q$ and $\rho=Q$. Theorem \ref{rhoto1mle4}, which is on $\rho\neq Q$, has its leading term in $\Lambda_I(\rho^k)$ involved with global information of the manifold, while the leading term in Theorem \ref{rhotome4}, which corresponds to $\rho=Q$, only depends on the geometric information at the blowup point. Moreover, the error terms in Theorem \ref{rhoto1mle4} and Theorem \ref{rhotome4} respectively also indicate the different asymptotic behaviors of blowup solutions near the singularity. All these differences in the two cases will lead to separate strategies in the construction of bubbling solutions in forthcoming works.

Since Liouville systems and Toda systems share a lot of common features, it is informative to compare
our main theorems with the ones for the $SU(3)$ Toda system. First, the location of the blowup point in Theorem \ref{locations} is a critical point of a combination of $\log h_i$, $\rho_i$ and $\gamma$ ($\nabla_1 \gamma$ vanishes if the Riemann surface has constant curvature). However for the $SU(3)$ Toda system, the blowup point $p$ is a critical point of both $\log h_1$ and $\log h_2$, i.e. $p$ satisfies (see \cite{lwz})
$$\nabla h_1(p)=\nabla h_2(p)=0. $$
Second, for the $SU(3)$ Toda system, the convergence rate of $\rho_i^k-\rho_i$ is estimated to be
$$\rho_i^k-\rho_i=(e^{\tilde c_i}b_i+o(1))\epsilon_k^2\log \epsilon_k^{-1}, $$
where $b_i$ is the term in (\ref{12apr23e10}). Nevertheless our result in (\ref{12apr23e10}) is again a combination of the $b_i$s. The comparison of the results reflects some major difference between the Toda system and our Liouville system:
\begin{enumerate}
\item
The dimension of kernel space of the linearized operator at an entire solution is $8$ for $SU(3)$ Toda system, and is $3$ for our Liouville system.
\item
The set $\Gamma:=\{(\rho_1,...,\rho_n);\quad
\rho_i=\int_{\mathbb R^2}e^{u_i},\quad (u_1,..,u_n) \mbox{ is an entire solution}\,\, \}$ is only a point for $SU(3)$ Toda system, while for the Liouville system it is a $(n-1)$ dimensional manifold.
\end{enumerate}

As far as the blowup analysis is concerned, our Liouville system has disadvantages in both respects, as the kernel space is too small and $\Gamma$ is too large. For a sequence of bubbling solutions, it is extremely difficult to pin-point suitable approximating solutions from $\Gamma$, because at the beginning, the local energy of bubbling solutions could be estimated in some rough way. This rough estimate of the local energy leads to a small perturbation of global bubbling solutions. This perturbation on global solutions, albeit small,
has a non-negligible effect on the approximation of blowup solutions. This difficulty
is particularly evident when we study bubbling solutions with multiple blowup points in \cite{linzhang4}. Therefore our method to obtain those sharp results is different from ones in Chen-Lin \cite{ChenLin1} for the mean field equations and Lin-Wei-Zhao \cite{lwz} for the $SU(3)$ Toda system (The methods in \cite{ChenLin1} and \cite{lwz} are similar).

The organization of the paper is as follows. In section two we first prove a uniqueness theorem for globally defined linearized Liouville systems. This result plays a central role for the delicate blowup analysis in sections three to five. The main idea of the proof uses a monotonicity property of solutions and we introduce a way to use maximum principles suitable for Liouville systems. In the second part of section two, we study the asymptotic behavior of global solutions to the Liouville system on $\mathbb R^2$ and obtain some Pohozaev identities. In section three and section four we obtain a sharp expansion result for blowup solutions around a blowup point. Then in section five, for local equations we use Pohozaev identity to determine the locations of blowup points. Then in section six we return to the equation on manifold and compute the leading term for $\rho^k\to \rho$ in both situations ($\rho\neq Q$ or $\rho=Q$) and
complete the proofs of the main theorems.

\medskip

\noindent{\bf Acknowledgement}
Part of the paper was finished when the second author was visiting Taida Institute for Mathematical Sciences (TIMS) in March 2011 and June 2012. He is very grateful to TIMS for their warm hospitality. He also would like to thank the National Science Foundation for the partial support (NSF-DMS-1027628).

\section{Preliminary results on the Liouville systems}

In this section we prove two theorems on the Liouville systems with the matrix $A$ satisfying $(H1)$. They are important for the blowup analysis and the computation of the leading terms of $\rho\to \Gamma_1$ in this paper and $\rho\to \Gamma_N$ in the forthcoming work \cite{linzhang4}.

\subsection{A uniqueness theorem for the linearized system of $n$ equations}
In the first subsection we prove a uniqueness theorem for the linearized system of $n$ equations.

\begin{thm}\label{uniqlin}
Let $A$ be a matrix that satisfies (H1), $u=(u_1,..,u_n)$ be a
radial solution of
$$\left\{\begin{array}{ll}
-\Delta u_i=\sum_{j=1}^na_{ij}e^{u_j}\,\, \mbox{in}\,\, \mathbb R^2,\quad i\in I=\{1,..,n\}\\
\int_{\mathbb R^2}e^{u_i}<\infty.
\end{array}
\right.
$$
Suppose $\phi=(\phi_1,..,\phi_n)$ satisfies
$$-\Delta \phi_i=\sum_{j=1}^na_{ij}e^{u_j}\phi_j,\quad i\in I.$$
\begin{enumerate}
\item
\begin{equation}\label{1027e1}
|\phi_i(x)|\le C(1+|x|)^{\tau},\quad x\in \mathbb R^2,
\end{equation}
for some $\tau\in (0,1)$ and
$$\phi_i(0)=0,\quad i\in I. $$
Then there exist $c_1,c_2\in \mathbb R$ such that
$$\phi_i(x)=c_1u_i'(r)\frac{x_1}r+c_2u_i'(r)\frac{x_2}r \quad \mbox{in}\quad \mathbb R^2, \quad i\in I. $$
\item
If $
|\phi_i(x)|\le C$ for all $x\in \mathbb R^2$, then there exist $c_0,c_1,c_2\in \mathbb R$ such that
$$\phi_i(x)=c_0(ru_i'(r)+2)+c_1u_i'(r)\frac{x_1}r+c_2u_i'(r)\frac{x_2}r,\quad \mathbb R^2 ,\quad i\in I. $$
\item If $\phi_i(x)=O(|x|^2)$ near $0$ and $|\phi_i(x)|\le C(1+|x|)^{2-\epsilon_0}$ for some $\epsilon_0>0$, then $\phi_i\equiv 0$.
\end{enumerate}
\end{thm}

Before the proof of Theorem \ref{uniqlin} we first establish
\begin{lem}\label{1129lem1} Let $A$ and $u$ be as in Theorem \ref{uniqlin}, let $\Phi=(\Phi_1,..,\Phi_n)$ solve
\begin{equation}\label{12e1}
\left\{\begin{array}{ll}
\Phi_i''(r)+\frac 1r\Phi_i'(r)-\frac 1{r^2}\Phi_i(r)+\sum_{j=1}^na_{ij}e^{u_j}\Phi_j=0,\quad 0<r<\infty, \\
\\
|\Phi_i(r)|\le Cr/(1+r)^{\epsilon_0} \quad \mbox{ for some }\epsilon_0\in (0,1),  \quad \forall i\in I.
\end{array}
\right.
\end{equation}
Then there exists a constant $C$ such that $\Phi_i(r)=Cu_i'(r)$ for $i\in I$.
\end{lem}

\noindent{\bf Proof of Lemma \ref{1129lem1}:}

The proof is in two steps. First we show that under the assumption of $\Phi_i$ at infinity we have the following sharper decay estimate:
\begin{equation}\label{phidecay}
|\Phi_i(r)|\le Cr(1+r)^{-2},\quad 0<r<\infty,\quad i\in I.
\end{equation}

Indeed, let $\bar \psi_i(r)=\Phi_i(r)/r$. By direct computation we see that $\bar \psi=(\bar \psi_1,...,\bar \psi_n)$ satisfies
\begin{equation}\label{12e2}
\bar \psi_i''(r)+\frac 3r\bar \psi_i'(r)+\sum_ja_{ij}e^{u_j}\bar \psi_j=0, \quad 0<r<\infty.
\end{equation}
Clearly in order to show (\ref{phidecay}) we only need to show
$|\bar \psi_i(r)|\le Cr^{-2}$ for $r>1$ under the assumption that $|\bar \psi_i(r)|\le Cr^{-\epsilon_0}$ for $i\in I$ and $r>1$.
Let
$$\tilde \psi_i(t)=\bar \psi_i(e^t)\quad \mbox{ and }\quad \tilde u_i(t)=u_i(e^t)+2t, $$
it is easy to see that $\tilde \psi(t)=(\tilde \psi_1(t),..,\tilde
\psi_n(t))$ satisfies
\begin{equation}\label{psidc1}
\tilde \psi_i''(t)+2\tilde \psi_i'(t)+\sum_{j\in I}a_{ij}e^{\tilde u_j(t)}\tilde \psi_j(t)=0,\quad -\infty<t<\infty
\end{equation}
and our goal is to show
\begin{equation}\label{psidc2}
\tilde \psi_i(t)=O(e^{-2t})
\end{equation}
 knowing $\tilde \psi_i(t)=O(e^{-\epsilon_0t})$   for  $t>0$, $i\in I$. Set
 $$l_i=\frac 1{2\pi}\sum_{j\in I}a_{ij}\int_{\mathbb R^2}e^{u_j},\quad l=\min\{l_1,..,l_n\}. $$
  By Lemma \ref{globalcpt} below $l>2$. Let
$$h_i(t)=-\sum_{j\in I}a_{ij}e^{\tilde u_j(t)}\tilde \psi_j(t)=O(e^{(2-l-\epsilon_0)t}),\quad t>0. $$
Then
$$\tilde \psi_i(t)=C_0+C_1 e^{-2t}+\frac 12\int_0^th_i(s)ds-\frac 12e^{-2t}\int_0^te^{2s}h_i(s)ds. $$
Using the asymptotic rate of $h_i(t)$ at infinity we further have
$$\tilde \psi_i(t)=(C_0+\frac 12\int_0^{\infty}h_i(s)ds)+C_1e^{-2t}+O(e^{(2-l-\epsilon_0)t}). $$
Since $\tilde \psi_i(t)$ tends to $0$ as $t$ goes to infinity we know $\tilde \psi_i(t)=O(e^{-2t})$ if
$l+\epsilon_0\ge 4$, in which case (\ref{psidc2}) is established. Otherwise we obtain $\tilde \psi_i(t)=O(e^{(2-l-\epsilon_0)t})$. In the latter case, we apply the same procedure to obtain a better decaying rate of $\tilde \psi_i(t)$ at infinity. After finite steps, (\ref{psidc2}) is established.

\medskip

In the second step we complete the proof of the Lemma \ref{1129lem1}.
By way of contradiction we suppose there is a solution $\Phi=(\Phi_1,..,\Phi_n)$ that satisfies (\ref{12e1}) and
$\Phi$ is not a multiple of $f=(u_1'(r),..,u_n'(r))$.
Let $\psi_i^0=-u_i'(r)/r$, then clearly
both $\psi^0=(\psi_1^0,..,\psi_n^0)$ and $\bar \psi=(\bar \psi_1,..,\bar \psi_n)$ satisfy (\ref{12e2}).
We verify by direct computation that
$$\int_0^re^{u_j}\psi_j^0(\sigma)\sigma^3ds=-\int_0^r(e^{u_j})'\sigma^2d\sigma
=-e^{u_j(r)}r^2+2\int_0^re^{u_j}\sigma d\sigma>0. $$
Note that the last inequality is justified by $u'_i(r)<0$ for $r>0$ and $i\in I$. Also, since $e^{u_i}\le Cr^{-2-\delta}$ for some
$\delta>0$ and $r>1$, $\displaystyle{\int_0^{\infty}e^{u_i}\psi_i^0(\sigma)\sigma^3d\sigma<\infty}$. Based on the computation above we set
$$S=\bigg \{ t|\, \psi_i^t=\psi_i^0+t(\psi_i^0-\bar \psi_i) \, \mbox{and }\,
\int_0^re^{u_i}\psi_i^t(\sigma )\sigma^3d\sigma>0,\,\, \forall r>0,\,\, i\in I\,  \bigg \}. $$
Let $\bar \psi=(\bar \psi_1,..,\bar \psi_n)$ and $\psi^0=(\psi_1^0,..,\psi_n^0)$. We first observe that $\psi_i^0(0)>0$. Suppose
$\bar \psi \not\equiv \psi^0$, we can assume that $\bar \psi_1(0)\neq 0$ and $|\bar \psi_i(0)|<\psi_i^0(0)$ for all $i\in I$.

 From the definition of $S$ we immediately see that $0\in S$. Moreover, since $|\bar \psi_i(r)|\le Cr^{-2}$ near infinity, we
can choose $|t|<\delta$ with $\delta$ small so that all $|t|<\delta$ belong to $S$. Another immediate observation is that $S$ has a lower
bound. Indeed, for $T$ sufficiently negative, $\psi_1^T(0)<0$, which is impossible to have $\int_0^re^{u_1}\psi_1^T(\sigma)\sigma^3d\sigma>0$ for $r$ small.

Let $T=\inf S$ and let $t_m\to T^+$. The sequence
$\psi^{t_m}_i$ obviously converges to a function $\psi_i$, which is just $\psi_i^0+T(\psi_i^0-\bar \psi_i)$. $\psi=(\psi_1,..,\psi_n)$ satisfies the property
\begin{equation}\label{12e4}
\int_0^re^{u_i(\sigma)}\psi_i(\sigma)\sigma^3d\sigma\ge 0,\quad \forall r>0,\quad i=1,..,n.
\end{equation}
On the other hand, from the behavior of $\bar \psi$ and $\psi^0$ at infinity (both are $O(r^{-2})$) we immediately observe that
$$\int_0^{\infty}e^{u_i(r)}\psi_i(r)r^3dr<\infty. $$
In regard to (\ref{12e2}) we have
\begin{equation}\label{12e3}
r^3\psi_i'(r)=-\int_0^r\sum_{j}a_{ij}e^{u_j}\psi_j\sigma^3 d\sigma, \quad 0<r<\infty.
\end{equation}
From (\ref{12e4}) and (\ref{12e3}) we see that $\psi_i$ is non-increasing. Since we have known that $|\psi_i(r)|\le Cr^{-2}$ near
infinity we have
\begin{equation}\label{12e6}
\psi_i(r)\ge 0, \quad \forall r>0, \quad i\in I.
\end{equation}

It is not possible to have all $\psi_i(0)=0$ because this implies $\psi_i\equiv 0$, a contradiction
to the assumption that $\bar \psi$ is not a multiple of $\psi^0$. Therefore without loss of generality we assume $\psi_1(0)>0$. Then we
further claim that $\psi_i$ is strictly decreasing for all $i\in I$. Indeed, let $I_1=\{ j\in I| \, a_{1j}>0 \, \}$, for each $j\in I_1$
we use (\ref{12e6}) and $\psi_1(0)>0$ to obtain
$$r^3\psi_j'(r)\le -\int_0^ra_{1j}e^{u_1}\psi_1\sigma^3d\sigma<0, \quad 0<r<\infty. $$
Therefore for each $j\in I_1$, $\psi_j$ is strictly decreasing, which immediately implies that $\psi_j(0)>0$. We can further define
$I_2=\{ i\in I| \, a_{ij}>0 \mbox{ for some }j\in I_1. \}$. Then the same argument shows that $\psi_i$ is strictly decreasing for each
$i\in I_2$ as well. Since the matrix $A=(a_{ij})_{n\times n}$ is irreducible, this process exhausts all $i\in I$.

(\ref{12e3}) yields
$\psi_i'(r)\le -Cr^{-3}$ for $r>1$ and $i\in I$. Then by
using $\lim_{r\to \infty}\psi_i(r)=0$ we further have
\begin{equation}\label{12e5}
\psi_i(r)\ge Cr^{-2},\quad r\ge 1.
\end{equation}
Then it is easy to see that for $t=T-\epsilon$ with $\epsilon>0$ small, we also have
$$\int_0^re^{u_i}\psi_i^t(\sigma)
\sigma^3d\sigma>0,\quad \mbox{for all } r>0, $$
a contradiction to the definition of $T$. Lemma \ref{1129lem1} is proved. $\Box $

\medskip

\noindent{\bf Proof of Theorem \ref{uniqlin}:} We first prove the third statement. The following function plays an important role:
 Let $f=(f_1,..,f_n)=(u_1',..,u_n')$. Then
\begin{equation}\label{85e1}
-\Delta f_i=\sum_ja_{ij}e^{u_j}f_j-\frac 1{r^2}f_i,\quad i=1,..,n.
\end{equation}

Let $\phi^k=(\phi_1^k,..,\phi_n^k)$ be defined as
$$\phi_i^k(r)=\frac 1{2\pi}\int_0^{2\pi}\phi_i(r\cos \theta,r\sin \theta)\cos k\theta d\theta,\quad i\in I. $$
Then $\phi^k$ satisfies
\begin{equation}\label{85e2}
-\Delta \phi_i^k=\sum_ja_{ij}e^{u_j}\phi_j^k-\frac{k^2}{r^2}\phi_i^k,\quad i\in I, \quad k=2,..
\end{equation}
Clearly $\phi_i^k(r)=o(r)$ near $0$ and $\phi_i^k(r)=O(r^{2-\epsilon_0})$ at $\infty$. We claim that
\begin{equation}\label{11mar29e1}
\phi_i^k(r)\equiv 0,\quad \forall k\ge 2,\quad \mbox{provided that }\,\, \phi_i^k(r)=O(r^{k-1+\tau}), \,\, r>1,\,\,  k\ge 2.
\end{equation}
Note that the growth condition in (\ref{11mar29e1}) is weaker than what is assumed in the assumption in Theorem \ref{uniqlin}.

The argument below also applies if $\phi$ is projected on $\sin k\theta$.
First we show that $\phi_i^k=o(r^{-1})$ as
$r\to \infty$. Indeed, using $\phi_j^k(x)\le C|x|^{k-1+\tau}$ we write $\sum_ja_{ij}e^{u_j}\phi^k_j$ as
$O(r^{k-1+\tau-2-\delta_0})$ (for some $\delta_0>0$). Let $g(t)=\phi_i^k(e^t)$, then from (\ref{85e2})
$g(t)$ satisfies
$$g''(t)-k^2g(t)=h(t),\quad t\in \mathbb R $$
where
\begin{equation}\label{11mar29e2}
h(t)=O(e^{(k-1+\tau-\delta_0)t})\quad  t>0.
\end{equation}
Let $g_1(t)=e^{kt}$ and $g_2(t)=e^{-kt}$ be two fundamental solutions of the homogeneous equation,
a general solution $g(t)$ is of the form
$$
g(t)=c_1g_1(t)+c_2g_2(t)-\frac{g_1(t)}{2k}\int_0^tg_2(s)h(s)ds+\frac {g_2(t)}{2k}\int_0^tg_1(s)h(s)ds
$$
where $c_1,c_2$ are constants.
Using (\ref{11mar29e2}) in the above we obtain
$$g(t)=c'_1g_1(t)+c'_2g_2(t)+O(e^{(k-1+\tau-\delta_0)t}),\quad \mbox{for }t>1$$
where $c'_1,c'_2$ are two constants.
Since $g(t)=O(e^{(k-1+\tau )t})$ for $t\to \infty$, we see that $c'_1=0$ and therefore $g(t)=O(e^{(k-1+\tau-\delta_0)t})$
as $t\to \infty$. Equivalently
\begin{equation}\label{11mar29e3}
\phi_i^k(r)=O(r^{k-1+\tau-\delta_0}) \quad i\in I.
\end{equation}
With (\ref{11mar29e3})
we further obtain
$$\sum_ja_{ij}e^{u_j}\phi^k_j=O(r^{k-1+\tau-2-2\delta_0}). $$ Consequently
$\phi_i^k(x)=O(r^{k-1+\tau-2\delta_0})$. Keep doing this for finite steps we obtain that $\phi_i^k$ decays faster than $r^{-1}$ at infinity. The asymptotic theory of ODE can be
similarly used to show that $\phi_i^k(r)=o(r)$ as $r\to 0$.

To get a contradiction, without loss of generality, we may assume that some of $\phi_i^k$, say $\phi_1^k(r)>0$ for some $r>0$ and
$$\max_{\mathbb R^+}\bigg (\frac{\phi_1^k(r)}{f_1(r)}\bigg )
=\max_{1\le j\le n}\bigg (\max_{\mathbb R^+}\frac{\phi_j^k(r)}{f_j(r)}\bigg ). $$
By noting $\phi_1^k(r)=o(r)$ as $r\to 0$ and $\phi_1^k(r)=o(\frac 1r)$ as $r\to \infty$, $\phi_1^k/
f_1(r)$ attains its maximum at some point $r_0\in \mathbb R^+$. Let $w_1(r)=\phi_1^k(r)/f_1(r)$. By a
direct computation, $w_1(r)$ satisfies
\begin{equation}\label{1019e1}
\Delta w_1+2\nabla w_1\cdot \frac{\nabla f_1}{f_1}+\frac{1-k^2}{r^2}w_1=\sum_{j=2}^na_{1j}e^{u_j}
\bigg (\frac{w_1f_j-\phi_j}{f_1}\bigg ).
\end{equation}
Now we apply the maximum principle at $r=r_0$, and obtain
$$\Delta w_1(r_0)\le 0,\quad \mbox{and}\quad \nabla w_1(r_0)=0. $$
Since $k>1$, (\ref{1019e1}) yields
\begin{equation}\label{1019e2}
\sum_{j=2}^na_{1j}e^{u_j}\bigg (\frac{w_1f_j-\phi_j^k}{f_1}\bigg )(r_0)<0
\end{equation}
because $w_1(r_0)>0$. On the other hand, for $j\ge 2$,
$$w_1(r_0)f_j(r_0)-\phi_j^k(r_0)=f_j(r_0)(\frac{\phi_1^k(r_0)}{f_1(r_0)}-\frac{\phi_j^k(r_0)}{f_j(r_0)})\ge
0,$$
which obviously contradicts (\ref{1019e2}). Therefore (\ref{11mar29e1}) is established. When $k=1$, $\phi_i^1\equiv 0$ because by
Lemma \ref{1129lem1}, $\phi_i^1(r)=Cu_i'(r)$. By the assumption $\phi_i^1(x)=O(|x|^2)$ near $0$, $C=0$.
The third statement of Theorem \ref{uniqlin} is established.

Again by Lemma \ref{1129lem1} the first statement of Theorem \ref{uniqlin} is established.

Finally, the second statement of Theorem \ref{uniqlin} is an immediate consequence of Lemma 3.1 of \cite{linzhang1}. Theorem \ref{uniqlin}
is established. $\Box$

\bigskip

\subsection{A Pohozaev identity for global solutions}

\begin{lem}\label{globalcpt}
Let $u=(u_1,...,u_n)$ be an entire, radial solution of
$$\left\{\begin{array}{ll}-\Delta u_i=\sum_{j=1}^na_{ij}e^{u_j}, \quad \mbox{in}\quad \mathbb R^2, \\
\int_{\mathbb R^2}e^{u_i}<\infty
\end{array}
\right.
$$
where $A$ is a constant matrix that satisfies $(H1)$. Let
$$c_i=u_i(0)+\frac 1{2\pi}\int_{\mathbb R^2}\log |\eta |(\sum_{j=1}^n a_{ij}e^{u_j(\eta)})d\eta, $$
$$\sigma_{i}=\frac 1{2\pi}\int_{\mathbb R^2}e^{u_i},\quad l_i=\sum_{j=1}^na_{ij}\sigma_{j},\,\, l=\min\{l_{1},...,l_{n}\} $$
and
$$\sigma_{iR}=\frac 1{2\pi}\int_{B_R}e^{u_i}. $$
Then for some $\delta_0>0$
\begin{equation}\label{12apr29e1}
e^{u_i(r)}=e^{c_i}r^{-l_i}(1+o(r^{-\delta_0})), \quad r>1,
\end{equation}
\begin{equation}\label{1105e1}
4\sum_{i\in I}\sigma_{iR}=\sum_{i,j\in I}a_{ij}\sigma_{iR}\sigma_{jR}+2\sum_{i\in I}e^{c_{i}}R^{2-l_{i}}+O(R^{2-l-\delta_0}).
\end{equation}
\end{lem}

\noindent{\bf Proof of Lemma \ref{globalcpt}:}

It is well known that
\begin{equation}\label{11may27e2}
u_i(x)=-\frac 1{2\pi}\int_{\mathbb R^2}\log |x-\eta |(\sum_{j=1}^n a_{ij}e^{u_j(\eta)})d\eta+c_i.
\end{equation}
Indeed, let $w_i$ be the function defined by the right hand side of (\ref{11may27e2}). Then $w_i-u_i$ is a harmonic function. Since they both have logarithmic growth at infinity, $w_i-u_i=c$. Evaluating both functions at $0$ we have
$c=c_i$.

Clearly
\begin{equation}\label{12apr24e1}
u_i(x)+l_{i}\log |x|=-\frac 1{2\pi}\int_{\mathbb R^2}\log \frac{|x-\eta |}{|x|}\sum_ja_{ij}e^{u_j(\eta)}d\eta+c_{i}.
\end{equation}
Using $\sum_j a_{ij}e^{u_j(r)}=O(r^{-2-\delta_0})$ for some $\delta_0>0$ and $r$ large, we obtain, by elementary estimates,
$$u_i(r)=-l_i \log r+c_i+o(r^{-\delta_0}), $$
which leads to
\begin{equation}\label{11may30e1}
\sigma_{i}=\sigma_{iR}+\frac{e^{c_{i}}}{l_{i}-2}R^{2-l_{i}}+O(R^{2-l_{i}-\delta}).
\end{equation}
We arrive at (\ref{1105e1}) by using (\ref{11may30e1}) in the Pohozaev identity for $\sigma$:
$$4\sum_{i} \sigma_{i}=\sum_{i,j}a_{ij}\sigma_{i}\sigma_{j}. $$
Lemma \ref{globalcpt} is established. $\Box$.

\section{First order estimates}\label{blowup}

Let $(h_1^k,...,h_n^k)$ be a family of positive, $C^3$ functions on $B_1$ with a uniform bound on their positivity and $C^3$ norm:
\begin{equation}\label{hikn}
\frac 1C\le h_i^k(x)\le C,\quad \|h_i^k\|_{C^3(B_1)}\le C, \quad x\in B_1, \quad i=1,..,n.
\end{equation}
In the next three sections we consider a sequence of locally defined, fully blown-up solutions $u^k=(u_1^k,...,u_n^k)$
and we shall derive their precise asymptotic behavior near their singularity and the precise location of their singularity.
Here we abuse the notation $u^k=(u_1^k,...,u_n^k)$ and it is independent of the one used in the introduction.

Specifically we assume that $u^k$ satisfies the following equation in $B_1$, the unit ball:
\begin{equation}\label{uik}
-\Delta u_i^k=\sum_{j=1}^na_{ij}h_j^ke^{u_j^k},\quad i=1,..,n,\quad x\in B_1,
\end{equation}
with $0$ being the only blowup point in $B_1$:
$$\max_Ku_i^k\le C(K), \quad \forall K\subset\subset \bar B_1\setminus \{0\}, \,\, \mbox{ and } \max_{B_1}u_i^k\to \infty,$$
with bounded oscillation on $\partial B_1$:
\begin{equation}\label{88e1}
|u_i^k(x)-u_i^k(y)|\le C,\forall x,y\in \partial B_1,\quad C\mbox{ independent of }k
\end{equation}
and uniformly bounded energy:
\begin{equation}\label{89e1}
\int_{B_1}e^{u_i^k}\le C,\quad C \mbox{ is independent of } k.
\end{equation}
Finally we assume that $u^k$ is a fully blown-up sequence, which means when re-scaled according its maximum, $\{u_k\}$ converges to a system of $n$ equations: Let
$u_1^k(0)=\max_{B_1}u_1^k$ and $\epsilon_k=e^{-\frac 12u_1^k(0)}$,
and
\begin{equation}\label{vikdef}
v_i^k(y)=u_i^k(\epsilon_ky)-u_1^k(0),\quad y\in \Omega_k:=B(0,\epsilon_k^{-1}).
\end{equation}
Then  $v^k=(v_1^k,...,v_n^k)$ converges in $C^2_{loc}(\mathbb R^2)$
to $v=(v_1,..,v_n)$, which satisfies
\begin{equation}\label{817e1}
\left\{\begin{array}{ll}-\Delta v_i=\sum_{j}a_{ij}\tilde h_j(0)e^{v_j},\quad \mathbb R^2,\quad i=1,..,n\\
\\
\int_{\mathbb R^2}e^{v_i}<\infty,\quad i=1,..,n,\quad v_1(0)=0,
\end{array}
\right.
\end{equation}
where $\tilde h_j(0)=\lim_{k\to \infty}h_j^k(0)$.

\bigskip

For the rest of the paper we set
$$m_i:=\frac{1}{2\pi}\int_{\mathbb R^2}\sum_{j=1}^na_{ij} \tilde h_j(0)e^{v_j}>2,\quad i\in I $$
and $m=\min\{m_1,...,m_n\}$.  In this section we derive a first order estimate for $v^k$ in $\Omega_k$.
In \cite{linzhang1} the authors prove that there is a sequence of global radial solutions $U^k=(U_1^k,..,U_n^k)$ of (\ref{817e1}) such that
 \begin{equation}\label{11mar9e1}
 |u_i^k(\epsilon_ky)-U_i^k(y)|\le C,\quad \mbox{for } |y|\le r_0\epsilon_k^{-1}.
 \end{equation}
 From (\ref{11mar9e1}) we have the following spherical Harnack inequality:
 \begin{equation}\label{11mar9e3}
 |u_i^k(\epsilon_ky)-u_i^k(\epsilon_k y')|\le C
 \end{equation}
 for all $|y|=|y'|=r\le r_0\epsilon_k^{-1}$ and $C$ is a constant independent of $k,r$. (\ref{11mar9e3}) will play an essential role in the first order estimate.
To improve (\ref{11mar9e1}) we introduce $\phi_i^k$ to be a harmonic function:
\begin{equation}\label{816e5}
\left\{\begin{array}{ll}
-\Delta \phi_i^k=0,\quad B_1,\\
\\
\phi_i^k=u_i^k-\frac 1{2\pi}\int_{\partial B_1}u_i^k,\quad \mbox{on}\quad \partial B_1.
\end{array}
\right.
\end{equation}
Obviously $\phi_i^k(0)=0$ by the mean value theorem and $\phi_i^k$ is uniformly bounded on $B_{1/2}$ because of (\ref{88e1}). Later in section 6, when the results in section 3,4,5 will be used to prove the main theorems, the function
$\phi_i^k$ will be specified when we consider the system defined on Riemann surface.

 Let $V^k=(V_1^k,..,V_n^k)$ be the radial solutions of
\begin{equation}\label{88e2}
\left\{\begin{array}{ll}-\Delta V_i^k=\sum_{j=1}^na_{ij}h_j^k(0)e^{V_j^k}\quad \mathbb R^2,\quad i\in I\\  \\
V_i^k(0)=v_i^k(0),\quad i\in I
\end{array}
\right.
\end{equation}
where $v_i^k$ is defined in (\ref{817e1}). It is easy to see that any radial solution $V$ of (\ref{88e2}) exists for all $r>0$ and $e^{V_i}\in L^2(\mathbb R^2)$.
The main result of this section is to prove that $V_i^k(y)+\phi_i^k(\epsilon_ky)$ is the first order approximation to $v_i^k(y)$.
 \begin{thm}\label{thm2}
Let $A$, $u^k$, $h^k=(h_1^k,..,h_n^k)$ ,$\phi_i^k$ and $v^k$ be described as above.
Then for any $\delta>0$, there exist $k_0(\delta)>1$ and $C$ independent of $k$ and $\delta$ such that for all $k\ge k_0$,
\begin{eqnarray}\label{11mar9e2}
&&|D^{\alpha}(v_i^k(y)-V_i^k(y)-\phi_i^k(\epsilon_ky))|\\
&\le &\left\{\begin{array}{ll}C\epsilon_k(1+|y|)^{3-m+\delta-|\alpha|},\quad m\le 3,\\
C\epsilon_k(1+|y|)^{\delta-l},\quad m>3,
\end{array}
\right.
\quad |y|<\epsilon_k^{-1}/2,\quad |\alpha|=0,1,2. \nonumber
\end{eqnarray}
\end{thm}

\begin{Def}\label{11def1}
$$\sigma_i^k=\frac 1{2\pi}\int_{\mathbb R^2} h_i^k(0)e^{V_i^k},\,\, m_i^k=\sum_{j=1}^n a_{ij}\sigma_j^k,\,\, m^k=\min \{m_1^k,..,m_n^k\}. $$
From Theorem \ref{thm2} it is easy to see that $\lim_{k\to\infty}m_i^k=m_i$. Thus $m_i^k\ge 2+\delta_0$ for some $\delta_0>0$ independent of $k$.
\end{Def}

%\begin{rem}\label{vik0}
%For $|\alpha |=0$ in (\ref{11mar9e2}) we use $v_i^k(0)-V_i^k(0)-\phi_i^k(0)=0$ to further obtain
%$$ |(v_i^k(y)-V_i^k(y)-\phi_i^k(\epsilon_ky))|\le \left\{\begin{array}{ll}
%C\epsilon_k|y|(1+|y|)^{2-m+\epsilon},\quad m\le 3, \\
%C\epsilon_k|y|(1+|y|)^{\epsilon-1},\quad m>3.
%\end{array}
%\right.
%$$
%\end{rem}

To prove Theorem \ref{thm2}, we have
\begin{equation}\label{1020e1}
-\Delta (v_i^k(y)-\phi_i^k(\epsilon_ky))=\sum_ja_{ij}H_j^k(\epsilon_ky)e^{v_j^k(y)-\phi_j^k(\epsilon_ky)},\,\, \mbox{ in } \Omega_k\,\, \big( \mbox{ see (\ref{vikdef})} \big )
\end{equation}
where
\begin{equation}\label{Hik}
H_i^k(\cdot)=h_i^k(\cdot)e^{\phi_i^k(\cdot)}.
\end{equation}
Since $\phi_i^k(0)=0$ we have $H_i^k(0)=h_i^k(0)$. Also, the definition of $\phi_i^k$
implies that $v_i^k-\phi_i^k(\epsilon_k\cdot)$ is a constant on $\partial \Omega_k$.

To estimate the error term $w_i^k=v_i^k-\phi_i^k(\epsilon_k \cdot )-V_i^k$. We find $w_i^k$ satisfies
\begin{equation}\label{wik01}
\left\{\begin{array}{ll}
\Delta w_i^k(y)+\sum_ja_{ij}H_j^k(\epsilon_ky)e^{\xi_j^k}w_j^k=-\sum_ja_{ij}(H_j^k(\epsilon_ky)-H_j^k(0))e^{V_j^k}, \\ \\
w_i^k(0)=0,\quad i\in I,\quad \nabla w_1^k(0)=O(\epsilon_k),
\end{array}
\right.
\end{equation}
where
$\xi_i^k$ is defined by
\begin{equation}\label{1211e1}
e^{\xi_i^k}=\int_0^1e^{tv_i^k+(1-t)V_i^k}dt.
\end{equation}

Since both $v^k$ and $V^k$ converge to $v$, $w_i^k=o(1)$ over any compact subset of $\mathbb R^2$.
The first estimate of $w_i^k$ is the following
\begin{lem}\label{1020lem1}
\begin{equation}\label{820e1}
w_i^k(y)=o(1)\log (1+|y|)+O(1),\quad \mbox{ for } y\in \Omega_k.
\end{equation}
\end{lem}

\medskip

\noindent{\bf Proof}: By (\ref{11mar9e3})
$$|v_i^k(y)-\bar v_i^k(|y|)|\le C, \quad \forall y\in \Omega_k $$
where $\bar v_i^k(r)$ is the average of $v_i^k$ on $\partial B_r$:
$$\bar v_i^k(r)=\frac{1}{2\pi r}\int_{\partial B_r}v_i^k. $$
 Thus we have
$e^{v_i^k(y)}=O(r^{-2-\delta_0})$
and $e^{V_i^k(y)}=O(r^{-2-\delta_0})$ where $r=|y|$ and $\delta_0>0$. Then
$$r(\bar w_i^k)'(r)=\frac 1{2\pi }\bigg (\int_{B_r}\sum_ja_{ij}H_j^k(\epsilon_k\cdot)e^{v_j^k}
-\int_{B_r}\sum_ja_{ij}h_j^k(0)e^{V_j^k}\bigg )$$

It is easy to use the decay rate of $e^{v_i^k}$, $e^{V_i^k}$ and the closeness between $v_i^k$ and $V_i^k$
to obtain
$$ r(\bar w_i^k)'(r)=o(1),\quad r\ge 1. $$
Hence $\bar w_i^k(r)=o(1)\log r$ and (\ref{820e1}) follows from this easily. Lemma \ref{1020lem1}
is established. $\Box$

\medskip

The following estimate is immediately implied by Lemma \ref{1020lem1}:
$$e^{\xi_i^k(y)}\le C(1+|y|)^{-m+o(1)}\quad \mbox{ for } y\in \Omega_k=B(0,\epsilon_k^{-1}).$$

Before we derive further estimate for $w_i^k$ we establish a useful estimate for the Green's function on $\Omega_k$ with respect to the
Dirichlet boundary condition:

\begin{lem}\label{greenk} Let $G(y,\eta)$ be the Green's function with respect to Dirichlet boundary condition on $\Omega_k$. For $y\in \Omega_k$,
let \begin{eqnarray*}
\Sigma_1&=&\{\eta \in \Omega_k;\quad |\eta |<|y|/2 \quad \}\\
\Sigma_2&=&\{\eta \in \Omega_k;\quad |y-\eta |<|y|/2 \quad \}\\
\Sigma_3&=&\Omega_k\setminus (\Sigma_1\cup \Sigma_2).
\end{eqnarray*}
Then in addition for $|y|>2$,
\begin{equation}\label{1020e5}
|G(y,\eta)-G(0,\eta)|\le \left\{\begin{array}{ll}
C(\log |y|+|\log |\eta ||),\quad \eta\in \Sigma_1,\\
C(\log |y|+|\log |y-\eta ||),\quad \eta\in \Sigma_2,\\
C|y|/|\eta |,\quad \eta \in \Sigma_3.
\end{array}
\right.
\end{equation}
\end{lem}

\noindent{\bf Proof:}
The expression for $G(y,\eta)$ is
$$G(y,\eta)=-\frac{1}{2\pi}\log |y-\eta |+\frac 1{2\pi}\log (\frac{|y|}{\epsilon_k^{-1}}|\frac{\epsilon_k^{-2}y}{|y|^2}-\eta |), \quad y,\eta\in \Omega_k.$$
In particular
$$G(0,\eta)=-\frac{1}{2\pi}\log |\eta |+\frac 1{2\pi}\log \epsilon_k^{-1},\quad \eta\in \Omega_k. $$
Therefore we write $G(y,\eta)-G(0,\eta)$ as
\begin{equation}\label{1020e3}
G(y,\eta)-G(0,\eta)=\frac 1{2\pi}\log \frac{|\eta |}{|y-\eta |}
+\frac 1{2\pi}\log |\frac{y}{|y|}-\frac{|y|\eta}{\epsilon_k^{-2}}|.
\end{equation}

The proof of (\ref{1020e5}) for $\eta\in \Sigma_1$ is obvious. For $\eta \in \Sigma_2$, (\ref{1020e5}) also obviously holds if either $|y|$ or $|\eta |$ is less than $\frac 78 \epsilon_k^{-1}$ because in this case
$$\bigg |\log \big | \frac{y}{|y|}-\epsilon_k^2|y|\eta \big |\bigg |\le C. $$
Consequently
\begin{eqnarray*}
|G(y,\eta)-G(0,\eta)|&\le &C(\log |\eta |+\big |\log |y-\eta |\big |+C)\\
&\le &C(\log |\eta  |+\big |\log |y-\eta |\big |).
\end{eqnarray*}
Therefore for $\eta\in \Sigma_2$ we only need to consider the case when $|y|,|\eta|>\frac 78\epsilon_k^{-1}$. In this case it is immediate to observe that $$\bigg |\log |\frac{y}{|y|}-\epsilon_k^2|y|\eta |\bigg |<C, \quad \mbox{ if }\,\, \angle (\frac{y}{|y}, \frac{\eta}{|\eta |})>\frac{\pi}8 $$
where $\angle(\cdot,\cdot)$ is the angle between two unit vectors. Thus for $\eta\in \Sigma_2$ we only consider the situation when $|y|,|\eta|>\frac 78\epsilon_k^{-1}$, $\angle(\frac{y}{|y|},\frac{\eta}{|\eta |})<\frac{\pi}8$. For this case we estimate $G(y,\eta)-G(0,\eta)$ as follows:
$$|G(y,\eta)-G(0,\eta)|\le |G(y,\eta)|+|G(0,\eta)|$$
$$|G(0,\eta)|\le C\log |y|.$$
\begin{eqnarray*}
|G(y,\eta)|&\le &\frac 1{2\pi}\big |\log |y-\eta |\big |+\frac 1{2\pi}\log \frac 87+\frac 1{2\pi}\big |\log |\frac{\epsilon_k^{-2}y}{|y|^2}-\eta |\big |\\
&\le & C(\log |y|+\big |\log |y-\eta |\big |)
\end{eqnarray*}
where the last inequality holds because
$$|y-\eta |\le |\frac{\epsilon_k^{-2}y}{|y|^2}-\eta |<C|y|, $$
which implies
$$\big |\log |\frac{\epsilon_k^{-2}y}{|y|^2}-\eta | \big |\le C(\big |\log |y-\eta |\big |+\log |y|). $$
The second case of (\ref{1020e5}) (when $\eta\in \Sigma_2$) is proved.

For $\eta\in \Sigma_3$, we first consider when $|\eta|>2|y|$. In this case
$$|\log \frac{|\eta |}{|\eta -y|}|=|\log |\frac{\eta}{|\eta |}-\frac{y}{|\eta |}||\le C\frac{|y|}{|\eta |}.$$
For the second term, since $\eta,y\in \Omega_k$ and $|\eta |>2|y|$, we have $|y||\eta |<\frac 12 \epsilon_k^{-2}$, consequently
$$ \bigg |\log |\frac{y}{|y|}-\frac{|y|\eta }{\epsilon_k^{-2}}| \bigg |\le C|y||\eta |\epsilon_k^2\le C\frac{|y|}{|\eta |}. $$
So (\ref{1020e5}) is proved in this case. Now we consider $\frac{|y|}2\le |\eta |\le 2|y|$ and $|\eta -y|\ge \frac{|y|}2$. For the first
term we have
$$\bigg |\log \frac{|\eta |}{|y-\eta |}\bigg |\le \log 4\le C\frac{|y|}{|\eta |}. $$
For the second term, we want to show
\begin{equation}\label{1027e5}
\bigg |\log |\frac{y}{|y|}-\frac{|y|\eta }{\epsilon_k^{-2}}| \bigg |\le C\le C\frac{|y|}{|\eta |}.
\end{equation}
If either $|y|\le \frac{15}{16}\epsilon_k^{-1}$ or $|\eta |\le \frac{15}{16}\epsilon_k^{-1}$ we have
$$|\frac{y}{|y|}-\frac{|y|\eta}{\epsilon_k^{-2}}|\ge \frac{1}{16},$$
therefore (\ref{1027e5}) obviously holds. For
$\frac{15}{16}\epsilon_k^{-1}<|y|, |\eta |\le \epsilon_k^{-1}$, using $|y-\eta |>\frac 12|y|$ we obtain easily
$$\bigg |\frac{y}{|y|}-\frac{|y|\eta}{\epsilon_k^{-2}}\bigg |\ge \frac 38. $$
Therefore (\ref{1020e5}) is proved in all cases. Lemma \ref{greenk} is established. $\Box$

\medskip

\noindent{\bf Proof of Theorem \ref{thm2}:} First we prove (\ref{11mar9e2}) for $\alpha=0$. We consider the case $m\le 3$, the proof for the case $m>3$ is similar.
By way of contradiction, we assume
$$\Lambda_k:=\max_{y\in\Omega_k}\frac{\max_{i\in I}|w_i^k(y)|}{\epsilon_k(1+|y|)^{3+\delta-m}}\to \infty. $$
Suppose $\Lambda_k$ is attained at $y_k\in \bar \Omega_k$ for some $i_0\in I$. We thus define
$$\bar w_i^k(y)=\frac{w_i^k(y)}{\Lambda_k\epsilon_k(1+|y_k|)^{3+\delta-m}}. $$
Here we require $\delta$ to be small so that $m-2-\delta>0$ (Thus $3-m+\delta<1$).
It follows from the definition of $\Lambda_k$ that for $y\in \Omega_k$
\begin{equation}\label{11apr20e1}
|\bar w_i^k(y)|=\frac{|w_i^k(y)}{\Lambda_k\delta_k(1+|y|)^{3+\delta-m}}\frac{(1+|y|)^{3+\delta-m}}{(1+|y_k|)^{3+\delta-m}}\le \frac{(1+|y|)^{3+\delta-m}}{(1+|y_k|)^{3+\delta-m}}.
\end{equation}
The equation for $\bar w_i^k$ is
\begin{equation}\label{11apr20e2}
-\Delta \bar w_i^k(y)=\sum_ja_{ij}h_j^k(0)e^{\xi_j^k}\bar w_j^k+o(1)\frac{(1+|y|)^{1-m}}{(1+|y_k|)^{3+\delta-m}},
\quad \Omega_k
\end{equation}
for $i\in I$. Here $\xi_i^k$ is given by (\ref{1211e1}). $\xi_i^k$ converges to $v_i$ in $C^2_{loc}(\mathbb R^2)$.  Besides, we also have $\bar w_i^k(0)=0$ for all $i$ and $\nabla \bar w_1^k(0)=o(1)$. If a subsequence of $y_k$ stays bounded, then along a subsequence $\bar w^k=(\bar w_1^k,..,\bar w_n^k)$ converges to $\bar w=(\bar w_1,...,\bar w_n)$ that satisfies
$$\left\{\begin{array}{ll}
-\Delta \bar w_i=\sum_ja_{ij}h_j(0)e^{v_j}\bar w_j,\quad \mathbb R^2, \quad i\in I,\\
\bar w_i(0)=0,\quad \nabla \bar w_1(0)=0,\quad |\bar w_i(y)|\le C(1+|y|)^{3+\delta-m},\quad y\in \mathbb R^2.
\end{array}
\right.
$$
Thanks to (1) of Theorem \ref{uniqlin}
$$\bar w_i(x)=c_1\frac{\partial v_i}{\partial x_1}+c_2\frac{\partial v_i}{\partial x_2}. $$
Since $\nabla \bar w_{1}(0)=0$ we have $c_1=c_2=0$, thus $\bar w_i\equiv 0$ for all $i$.
On the other hand, the fact that $\bar w_{i_0}^k(y_k)=\pm 1$ for some $i_0\in I$
implies that $\bar w_{i_0}(\bar y)=\pm 1$ where $\bar y$ is the limit of $y_k$. This contradiction means that $y_k\to \infty$.
Next we shall show a contradiction if $|y_k|\to \infty$. By the Green's representation formula for $\bar w_i^k$,
$$\bar w_i^k(y)=\int_{\Omega_k}G(y,\eta)(-\Delta \bar w_i^k(\eta))d\eta + \bar w_i^k|_{\partial \Omega_k}$$
where $\bar w_i^k|_{\partial \Omega_k}$ is the boundary value of $\bar w_i$ on $\partial \Omega_k$ ( which is a constant). From (\ref{11apr20e1}) and (\ref{11apr20e2}) we have
$$|-\Delta \bar w_i^k(\eta)|\le \frac{C(1+|\eta |)^{3+\delta-2m}}{(1+|y_k|)^{3+\delta-m}}+
\frac{C(1+|\eta |)^{1-m+\delta}}{\Lambda_k(1+|y_k|)^{3+\delta-m}}. $$

Thus for some $i\in I$ we have
\begin{eqnarray} \label{11apr20e3}
&&1=|\bar w_i^k(y_k)-\bar w_i^k(0)|\\
&\le &C\int_{\Omega_k}|G(y_k,\eta)-G(0,\eta)|\bigg (\frac{(1+|\eta |)^{3+\delta-2m}}{(1+|y_k|)^{3+\delta-m}}
+\frac{(1+|\eta |)^{1-m+\delta}}{\Lambda_k(1+|y_k|)^{3+\delta-m}}\bigg ),\nonumber
\end{eqnarray}
where the constant on the boundary is canceled out.
To compute the right hand side of the above, we decompose the $\Omega_k$ as $\Omega_k=\Sigma_1\cup \Sigma_2\cup \Sigma_3$ as in Lemma \ref{greenk}. Using (\ref{1020e5}) we have
$$\int_{\Sigma_1\cup \Sigma_2}|G(y_k,\eta)-G(0,\eta)|(1+|\eta |)^{3+\delta-2m}d\eta=O(1)(\log |y_k|)(1+|y_k|)^{(5+\delta-2m)_+} $$
where
$$(1+|y_k|)^{\alpha_+}=\left\{\begin{array}{ll}(1+|y_k|)^{\alpha}, &\quad \alpha>0, \\
\log (1+|y_k|), &\quad \alpha=0, \\
1, &\quad \alpha<0.
\end{array}
\right.
$$
$$\int_{\Sigma_3}|G(y_k,\eta)-G(0,\eta)|(1+|\eta |)^{3+\delta-2m}d\eta=O(1)(1+|y_k|)^{5+\delta-2m}. $$
Hence
$$\int_{\Omega_k}|G(y_k, \eta)-G(0,\eta)|\frac{(1+|\eta |)^{3+\delta-2m}}{(1+|y_k|)^{3+\delta-m}}d\eta
=O(1)(1+|y_k|)^{2-m}. $$
Similarly we can compute the other term:
$$\int_{\Omega_k}|G(y_k, \eta)-G(0, \eta)|\frac{(1+|\eta |)^{1-m+\delta}}{\Lambda_k(1+|y_k|)^{3+\delta-m}}d\eta=O(1)\Lambda_k^{-1}(\log (1+|y_k|))^{-\frac{\delta}2}. $$
By the computations above we see that the right hand side of (\ref{11apr20e3}) is $o(1)$, a contradiction to the left hand side of (\ref{11apr20e3}). Thus (\ref{11mar9e2}) is established for $\alpha=0$. The estimates for $|\alpha|=1$ and $2$ follow
easily by scaling and standard elliptic estimates. Therefore Theorem \ref{thm2} is completely proved. $\Box$

\section{Second order estimates}

In this section we improve the estimates in Theorem \ref{thm2} for $m<4$ and $m=4$, respectively.
Let $p_{i,k}$ be the maximum point of $v_i^k(\cdot)-\phi_i^k(\epsilon_k\cdot)$. The following lemma estimates the location of $p_{i,k}$.
\begin{lem}\label{pik}
$p_{i,k}=O(\epsilon_k),\quad i\in I. $
\end{lem}

\noindent{\bf Proof :}
Applying Theorem \ref{thm2} to $v_i^k-\phi_i^k(\delta \cdot )$ on $B_1$:
\begin{eqnarray}\label{11mar12e2}
&&D^{\alpha}(v_i^k(y)-\phi_i^k(\epsilon_ky))\\
&=&D^{\alpha}(V_i^k(|y|))+O(\epsilon_k),\quad |y|<1,\quad |\alpha |=0,1,2. \nonumber
\end{eqnarray}
The equation for $V_i^k$ is
$$(V_i^k)''(r)+\frac 1r(V_i^k)'(r)+\sum_{j=1}^na_{ij}H_j^k(0)e^{V_j^k(r)}=0,\quad r>0. $$
From $(V_i^k)'(0)=0$, we see $\lim_{r\to 0}(V_i^k)'(r)/r=(V_i^k)''(0)$. Thus
\begin{equation}\label{11apr4e2}
(V_i^k)''(0)=-\frac 12\sum_{j=1}^na_{ij}H_j^k(0)e^{V_j^k(0)}<-C
\end{equation}
for some $C>0$ independent of $k$.
Since $p_{i,k}$ is the maximum point of $v_i^k(\cdot)-\phi_i^k(\epsilon_k\cdot)$, we deduce from (\ref{11mar12e2})
that
$(V_i^k)'(|p_{i,k}|)=O(\epsilon_k)$, thus from (\ref{11apr4e2}) we have $p_{i,k}=O(\epsilon_k)$.
Lemma \ref{pik} is established. $\Box$

\medskip

The main result in this section is to find the $\epsilon_k$ approximation to $v_i^k(\cdot)-\phi_i^k(\epsilon_k\cdot)$. It is most convenient to write the expansion around one of the $p_{i,k}$s. We choose $p_{1,k}$ and shall use $\Phi^k=(\Phi_1^k,..,\Phi_n^k)$ to denote the projection of $v_i^k(\cdot)-\phi_i^k(\epsilon_k\cdot)$ onto $span\{\sin \theta, \cos\theta\}$. i.e.
\begin{equation}\label{11aug11e1}
\Phi^k_i(r\cos\theta,r\sin\theta)=\epsilon_k(G_{1,i}^k(r)\cos \theta+G_{2,i}^k(r)\sin\theta), \quad i\in I
\end{equation}
with $G_{t,i}^k(r)$ ($t=1,2$) satisfying some ordinary differential equations to be specified later.

Set $v^{1,k}=(v^{1,k}_1,..,v^{1,k}_n)$ as
\begin{equation}\label{vi0kd}
v^{1,k}_i(\cdot)=v_i^k(\cdot+p_{1,k})-\phi_i^k(\epsilon_k(\cdot+p_{1,k}))
\end{equation}
in
\begin{equation}\label{omegai0kd}
\Omega_{1,k}:=\{\eta;\,\, \eta+p_{1,k}\in \Omega_k\,\, \}.
\end{equation}
Using $\nabla v_i^k(0)=O(\epsilon_k)$ (by Theorem \ref{thm2}) and $\phi_i^k(0)=0$  we observe that
\begin{eqnarray}\label{11mar12e1}
&&v^{1,k}_i(0)=v_i^k(p_{1,k})-\phi_i^k(\epsilon_k p_{1,k})\\
&=&v_i^k(0)+\nabla v_i^k(0)\cdot p_{1,k}+O(|p_{1,k}|^2)+O(p_{1,k}\epsilon_k)\nonumber\\
&=&v_i^k(0)+O(\epsilon_k^2) \nonumber
\end{eqnarray}
The equation that $v^{1,k}$ satisfies is (combining (\ref{11mar9e2}) and (\ref{1020e1}))
\begin{equation}\label{vpke}
\left\{\begin{array}{ll}
\Delta v_i^{1,k}+\sum_{j=1}^na_{ij}H_j^{1,k}(y)e^{v_j^{1,k}}=0,\quad \mbox{in}\,\, \Omega_{1,k}\\
\\
\nabla v_1^{1,k}(0)=0, \quad \nabla v_i^{1,k}(0)=O(\epsilon_k),\quad i=2,...,n
\end{array}
\right.
\end{equation}
where $H^{1,k}=(H^{1,k}_1,...,H^{1,k}_n)$ is defined by (see (\ref{Hik}))
\begin{equation}\label{hi0k}
H_i^{1,k}(\cdot)=H_i^k(\epsilon_k\cdot+\epsilon_kp_{1,k})=h_i^k(\epsilon_k\cdot +\epsilon_kp_{1,k})e^{\phi_i^k(\epsilon_k\cdot +\epsilon_kp_{1,k})}.
\end{equation}
Trivially
\begin{equation}\label{12feb16e1}
H_i^{1,k}(0)=h_i^k(0)+O(\epsilon_k^2).
\end{equation}

In the coordinate around $p_{1,k}$ and we seek to approximate $v_i^{1,k}$ in
$\Omega_{1,k}$.
The first term in the approximation of  $v^{1,k}$ is $V^k$. Here we note that the domain $\Omega_{1,k}$ is shifted from the ball $\Omega_k$ by $p_{1,k}$.

We shall use five steps to establish an approximation of $v_i^{1,k}$ without distinguishing $m=4$ or not.

\medskip

\noindent{\bf Step one:}

Let $w^{1,k}=(w^{1,k}_1,..,w^{1,k}_n)$ be the difference between $v^k$ and $V^k$:
$$w^{1,k}_i(y)=v^{1,k}_i(y)-V^{k}_i(|y|), \quad y\in \Omega_{1k}. $$
Taking the difference between (\ref{vpke}) and (\ref{88e2}), we have
$$
\Delta w_i^{1,k}+\sum_j a_{ij}H_j^{1,k}(y)e^{V_j^k+w_j^{1,k}}-\sum_j a_{ij}h_j^k(0)e^{V_j^k}=0, $$
which is
$$\Delta w_i^{1,k}+\sum_j a_{ij}h_j^k(0)e^{V_j^k}(\frac{e^{w_j^{1,k}}H_j^{1,k}(y)}{h_j^k(0)}-1)=0. $$
Here we observe that the oscillation of $V_i^k $ on $\partial \Omega_k$ is $O(\epsilon_k^2)$. Indeed, recall that $\Omega_{1,k}$ is the shift of the large ball $\Omega_k$ by $p_{1,k}$.
 Let $y_1,y_2\in
\partial \Omega_{1,k}$, one can find $y_3$ such that $|y_3|=|y_2|$ and $|y_3-y_1|\le C\epsilon_k$. Since
$(V_i^{k})'(r)\sim r^{-1}$ for $r>1$ and $|y_1|\sim \epsilon_k^{-1}$, we have
\begin{equation}\label{12jun19e1}
V_i^{,k}(y_1)-V_i^{k}(y_2)=V_i^{k}(y_1)-V_i^{k}(y_3)=O(\epsilon_k^2).
\end{equation}

With (\ref{12jun19e1}) we further write the equation for $w_i^{1,k}$ as
\begin{equation}\label{newwk}
\left\{\begin{array}{ll}
\Delta w_i^{1,k}+ \sum_ja_{ij}h_j^k(0)e^{V_j^{k}}w_j^{1,k} =E_i^k,\quad \Omega_{1,k}.\\
\\
w_i^{1,k}(0)=O(\epsilon_k^2),\,\, i\in I,\quad \bar w_i^{1,k}=O(\epsilon_k^2) \mbox{ on } \partial \Omega_{1,k},\\
\\
\nabla w_{1}^{1,k}(0)=0,\quad \nabla w_i^{1,k}(0)=O(\epsilon_k),\quad i\in I\setminus \{1\}
\end{array}
\right.
\end{equation}
where
\begin{equation}\label{ewik}
E_i^k=-\sum_ja_{ij}h_j^k(0)e^{V_j^{k}}\bigg (\frac{H_j^{1,k}(y)}{h_j^k(0)}e^{w_j^{1,k}}-1-w_j^{1,k}\bigg ).
\end{equation}
Similar to Theorem \ref{thm2} we also have
\begin{lem}\label{wi0ki}
For any $\delta>0$, there exists $k_0(\delta)>1$ such that for some $C>0$ independent of $k$ and $\delta$, the following estimate holds for all $k\ge k_0$:
\begin{equation}\label{wi0k}
|w^{1,k}_i(y)|\le \left\{\begin{array}{ll}
C\epsilon_k(1+|y|)^{3-m+\delta},\,\, m\le 3,\\
C\epsilon_k(1+|y|)^{\delta},\,\, m> 3,
\end{array}
\right.
 \quad y\in \Omega_{1,k}.
\end{equation}
\end{lem}

\noindent{\bf Proof:}
Using the definition of $v^{1,k}$ and Theorem \ref{thm2} we have
$$|v_i^{1,k}(y)-V_i^k(y+p_{1,k})|\le \left\{\begin{array}{ll}
C\epsilon_k (1+|y|)^{3-m+\delta},\quad m\le 3, \\
C\epsilon_k (1+|y|)^{\delta},\quad m>3.
\end{array}
\right.
$$
On the other hand we clearly have
$$|V_i^k(y)-V_i^k(y+p_{1,k})|\le C\epsilon_k(1+|y|)^{-1} $$
by mean value theorem and the estimate of $\nabla V_i^k$.
Lemma \ref{wi0ki} is established.
$\Box$

\medskip

Using Lemma \ref{wi0ki} and (\ref{12feb16e1}) we now rewrite $E_i^k$, clearly
$$E_i^k=-\sum_ja_{ij}e^{V_j^{k}}(H_j^{1,k}(y)-h_j^k(0))+(H_j^k(y)-h_j^k(0))w_j^{1,k}+O((w_j^{1,k})^2). $$
By Lemma \ref{wi0ki} and (\ref{12feb16e1}), the last two terms are $O(\epsilon_k^2(1+|y|)^{2-m})$ regardless whether $m\ge 3$ or not.
Thus
\begin{eqnarray}\label{ewik1}
E_i^k&=& -\sum_ja_{ij}e^{V_j^{k}}(H_j^{1,k}(y)-h_j^k(0))+O(\epsilon_k^2)(1+|y|)^{2-m}\\
&=&-\sum_ja_{ij}e^{V_j^{k}}(H_j^{1,k}(y)-H_j^{1,k}(0))+O(\epsilon_k^2)(1+|y|)^{2-m} \nonumber
\end{eqnarray}
where in the last step we used (\ref{12feb16e1}) again.

\medskip

\noindent{\bf Step Two: Estimate of the radial part of $w^{1,k}$:}

Let $g^{k,0}=(g_1^{k,0},...,g^{k,0}_n)$ be the radial part of $w^{1,k}$:
$$g_i^{k,0}(r)=\frac 1{2\pi}\int_0^{2\pi}w_i^{1,k}(r\cos \theta, r\sin \theta)d\theta. $$
Due to the radial symmetry of $V_i^k$,
$g^{k,0}$ satisfies
\begin{equation}\label{11apr8e1}
\left\{\begin{array}{ll}
L_i g^{k,0}
=-\frac{\epsilon_k^2}4\sum_ja_{ij}\Delta H_j^{1,k}(0)r^2e^{V_j^{k}}
+O(\epsilon_k^2)(1+r)^{\delta-m}\\
\\
g_i^{k,0}(0)=O(\epsilon_k^2),\quad i\in I, \quad \frac{d}{dr}g_{1}^{k,0}(0)=0.
\end{array}
\right.
\end{equation}
where ( for simplicity we omit $k$ in $L_i$)
$$L_ig^{k,0}=\frac{d^2}{dr^2}g_i^{k,0}+\frac 1r\frac{d}{dr}g_i^{k,0}+\sum_ja_{ij}h_j^{k}(0)e^{V^{k}_j}g_j^{k,0}$$
We claim that for $\delta>0$, there exists $k_0(\delta)>1$ such that for all $k\ge k_0$
\begin{equation}\label{11mar29e4}
|g_i^{k,0}(r)|\le C\epsilon_k^2(1+r)^{4-m+\delta},\quad 0<r<\epsilon_k^{-1}
\end{equation}
holds for some $C$ independent of $k$ and $\delta$.
So, $g_i^{k,0}(r)$ can be discarded as an error term.

To prove (\ref{11mar29e4}), we first observe that $m\le 4$ and by (\ref{11apr8e1})
$$|L_ig^{k,0}|\le C\epsilon_k^2(1+r)^{2-m+\delta/2}. $$
Let $f^k=(f_1^k,..,f_n^k)$ be the solution of
$$\left\{\begin{array}{ll}
\frac{d^2}{dr^2}f_i^k+\frac 1r\frac{d}{dr}f_i^k=L_ig^{k,0}, \\
\\
f_i^k(0)=\frac{d}{dr}f_i^k(0)=0.
\end{array}
\right.
$$
Then elementary estimate shows
\begin{equation}\label{12feb13e1}
|f_i^k(r)|\le C\epsilon_k^2(1+r)^{4-m+\delta}.
\end{equation}
If $m>\frac 52$, we claim
\begin{equation}\label{11apr5e1}
|g_i^{k,0}(r)-f_i^k(r)|\le C\epsilon_k^2(1+r)^{\delta}, \quad 0<r<\epsilon_k^{-1}.
\end{equation}

Indeed, let $\bar g^{k}=g^{k,0}-f^{k}$,
then clearly
$$\left\{\begin{array}{ll}
L_i\bar g^k=F_i^k, \\
\\
\bar g_i^k(0)=O(\epsilon_k^2), \quad \frac{d}{dr}\bar g_{1}^k(0)=0,
\end{array}
\right.
$$
where
$$F_i^k:=-\sum_ja_{ij}h_j^{k}(0)e^{V^{k}_j}f_j^k=O(\epsilon_k^2)(1+r)^{4-2m+\delta}. $$

By considering $\frac{d^2}{dr^2}+\frac 1r\frac{d}{dr}$ as $\Delta $ in $\mathbb R^2$ and $\bar g^k$ a solution with boundary
oscillation $0$ in $B(0,\epsilon_k^{-1})$, we obtain (\ref{11apr5e1}) by the argument for Theorem \ref{thm2}. Here we note that to apply Theorem \ref{uniqlin}, it is essential to require $6-2m+\delta<1$ (it holds if $m>5/2$), the estimate on the Green's function in Lemma \ref{greenk} and the condition $\frac{d}{dr}\bar g_{1}^k(0)=0$. Since $m\le 4$, $O(\epsilon_k^2)(1+r)^{\delta}$ is part of the error.

If $m\le \frac 52$ we apply the same ideas by adding more correction functions to $g^{k,0}$: Let $N$ satisfy $2+(2-m)N<1$,
we add $N$ correcting functions to make the right hand side of the equation of the order $O(\epsilon_k^2(1+r)^{(2-m)N+\delta})$. Note that each correction can be discarded as an error in the sense that they are smaller than the right hand side of (\ref{11jun14e5}). Using $2+2N-mN<1$ and the argument in the proof of Theorem \ref{thm2} we obtain (\ref{11mar29e4}).

\medskip

\noindent{\bf Step Three: The projection on $\sin\theta$ and $\cos \theta$}

In this step we consider the projection of $w^{1,k}$ over $\cos\theta$ and $\sin \theta$, respectively:
$$\epsilon_kG_{1,i}^k(r)=\frac 1{2\pi}\int_0^{2\pi}w_i^{1,k}(r,\theta)\cos\theta d\theta, \quad
\epsilon_kG_{2,i}^k(r)=\frac 1{2\pi}\int_0^{2\pi}w_i^{1,k}(r,\theta)\sin\theta d\theta. $$
Let
\begin{equation}\label{1224e2}
\Phi_i^k=\epsilon_kG_{1,i}^k(r)\cos \theta +\epsilon_kG_{2,i}^k(r) \sin \theta
\end{equation}
clearly $G_{1,i}^k$ and $G_{2,i}^k$ solve the following linear systems: For $0<r<\epsilon_k^{-1}$ and $t=1,2$
\begin{eqnarray}\label{1030e1}
(\frac{d^2}{dr^2}+\frac 1r\frac{d}{dr}-\frac{1}{r^2})G^k_{t,i}+\sum_j a_{ij}h_j^k(0)e^{V_j^k}G^k_{t,j}\\
=-\sum_ja_{ij}\partial_tH_j^{1,k}(0)re^{V_j^{k}}
+O(\epsilon_k)(1+r)^{2-m}. \nonumber
\end{eqnarray}

Then $\Phi^k$ solves
\begin{eqnarray}
&&\Delta \Phi_i^k+\sum_ja_{ij}h_j^k(0)e^{V_j^{k}}\Phi_j^k\\
\label{1224e3}
&=&-\epsilon_k\sum_ja_{ij}(\partial_1H_j^{1,k}(0)y_1+\partial_2H_j^{1,k}(0)y_2)e^{V_j^{k}}
+O(\epsilon_k^2)(1+|y|)^{2-m}. \nonumber
\end{eqnarray}

By the long behavior of $w^{1,k}$ (Lemma \ref{wi0ki}) we have
\begin{equation}\label{11mar21e1}
|G_{1,i}^k(r)|+|G_{2,i}^k(r)|\le \left\{\begin{array}{ll}C(1+r)^{3-m+\delta}, \,\, m\le 3, \\
C(1+r)^{\delta}, \,\, m>3.
\end{array}
\right.
\end{equation}
Note that $G_{t,i}^k(0)=0$ and $G_{t,i}^k(r)=O(r)$ near $0$.

\medskip

\noindent{\bf Step Four: Projection of $w^{1,k}$ onto higher frequencies }

Let
$$g^{k,l}=(g^{k,l}_1,...,g^{k,l}_n)$$
be the projection of $w^{1,k}$ on $\sin l\theta$. In this step we first establish a preliminary estimate for all these projections:

\begin{lem}\label{11jun14lem3}
There exist $l_0\ge 3$ and $C>0$
independent of $k,l$ such that
$$|g_i^{k,l}(r)|\le C\epsilon_k^2r^2,\quad 0<r<\epsilon_k^{-1}, \quad \forall l\ge l_0.  $$
\end{lem}

\noindent{\bf Proof:}

By (\ref{newwk}), (\ref{ewik1}) and Lemma \ref{wi0ki} $g^{k,l}$ satisfies
\begin{equation}\label{11mar31e1}
\left\{\begin{array}{ll}
(\frac{d^2}{dr^2}+\frac{1}r\frac{d}{dr}-\frac{l^2}{r^2})g_i^{k,l}+\sum_ja_{ij}h_j^k(0)e^{V_j^{k}}g_j^{k,l}=h_{ikl}\\
\\
g_i^{k,l}(0)=0.
\end{array}
\right.
\end{equation}
where, applying the Taylor expansion of $H_j^{i,k}$ up to the second order,
\begin{equation}\label{12mar9e1}
|h_{ikl}(r)|\le C\epsilon_k^2(1+r)^{2-m+\delta/2}
\end{equation}
for some $C$ independent of $k$ and $l$.
Thus
\begin{equation}\label{11apr5e9}
\Delta g_i^{k,l}-\frac{l^2}{r^2}g_i^{k,l}+\sum_ja_{ij}h_j^{k}(0)e^{V^{k}_j}g_j^{k,l}>-c_0\epsilon_k^2(1+r)^{2-m+\delta}
\end{equation}
for some $c_0>0$.
Let
$$g(r)=\frac{r^2}4\int_r^{\infty}\frac{c_0}s(1+s)^{2-m+\delta}ds+c_0\frac{r^{-2}}4\int_0^rs^3(1+s)^{2-m+\delta}ds+r^2. $$
Then clearly $g(r)>0$ for $r>0$, $g$ solves
\begin{equation}\label{11apr11e11}
\left\{\begin{array}{ll}
g''+\frac 1r g'-\frac{4}{r^2}g(r)=-c_0r^2(1+r)^{-m+\delta}, \quad r>0, \\
\\
g(0)=g'(0)=0.
\end{array}
\right.
\end{equation}
and
$$g(r)=r^2+O(r^{4-m+\delta}), \quad r>1, \quad g(r)\le Cr^2\log (\frac 1r+1),\quad r\le 1. $$
Clearly by Lemma \ref{wi0ki}
\begin{equation}\label{11apr5e7}
\epsilon_k^2g(\epsilon_k^{-1})>\max_{y\in \partial \Omega_{1,k}}|w^{1,k}_i(y)|\ge |g_i^{k,l}(\epsilon_k^{-1})|.
\end{equation}
The reason that we include $r^2$ in the definition of $g(\cdot)$ is because by Lemma \ref{wi0ki} we only know
$w_i^{1,k}(x)=O(\epsilon_k (1+|x|)^{\delta})$ for $m>3$.
Let $g^k=(g_1^k,..,g_n^k)=\epsilon_k^2(g,..,g)$, then it is easy to see that for $l_0$ sufficiently large and $l>l_0$
\begin{eqnarray}\label{11apr5e8}
&&\Delta g_i^k-\frac{l^2}{r^2}g_i^k+\sum_ja_{ij}h_j^{k}(0)e^{V^{k}_j}g_j^k\\
&=&-\epsilon_k^2c_0(1+r)^{2-m+\delta}+\sum_ja_{ij}h_j^{k}(0)e^{V^{k}_j}g_j^k-\frac{l^2-4}{r^2}g_i^k \nonumber\\
&\le &-\epsilon_k^2c_0(1+r)^{2-m+\delta}.\nonumber
\end{eqnarray}
To prove Lemma \ref{11jun14lem3} it is enough to show
\begin{equation}\label{11jun15e1}
|g_i^{k,l}(r)|\le C g_i^k, \quad 0<r<\epsilon_k^{-1},\quad l\ge l_0.
\end{equation}
with $C$ independent of $k$ and $l$.
To this end, we shall use (\ref{11apr5e9}) and the initial value
\begin{equation}\label{11apr5e10}
g_i^{k,l}(r)=o(r) \mbox{ near }0.
\end{equation}
In the following we use the argument in the proof of Theorem \ref{uniqlin}.
Let $\psi_i^k=g_i^k-g_i^{k,l}$, by (\ref{11apr5e7}),(\ref{11apr5e8}),(\ref{11apr5e9}) and (\ref{11apr5e10})
$$\left\{\begin{array}{ll}
\Delta \psi_i^k-\frac{l^2}{r^2}\psi_i^k+\sum_ja_{ij}h_j^{k}(0)e^{V^{k}_j}\psi_j^k<0,\quad 0<r<\epsilon_k^{-1}, \\
\\
\psi_i^k(r)=o(r)\mbox{ near }0, \quad \psi_i^k(\epsilon_k^{-1})>0 \mbox{ for all } i.
\end{array}
\right.
$$
Our goal is to prove $\psi_i^k(r)>0$ for all $0<r<\epsilon_k^{-1}$. To do this, let $f^k=(f^k_1,...,f^k_n)$ be the positive solution to the homogeneous
system:
$$\Delta f_i^k-\frac{1}{r^2}f_i^k+\sum_ja_{ij}h_j^{k}(0)e^{V^{k}_j}f_j^k=0, \quad r>0 $$
such that
$$f_i^k(r)\sim r, \mbox{ near } 0, \quad f^k_i(r)\sim \frac 1r \mbox{ near } \infty, \quad f_i^k(r)>0 \mbox{ for all } r>0. $$
If $\psi_i^k$ is not always non-negative we assume
$$\min_{\mathbb R^+}\frac{\psi_1^k(r)}{f_1^k(r)}=\min_{1\le i\le n}\min_{\mathbb R^+}\frac{\psi_i^k(r)}{f_i^k(r)}<0. $$
Suppose the minimum is attained at $r_k$. By the behavior of $\psi_i^k$ and $f_i^k$, we have $0<r_k<\epsilon_k^{-1}$. Let $w_1^k=\psi_1^k/f_1^k$,
then $w_1^k$ satisfies
$$\Delta w_1^k+2\nabla w_1^k\cdot \frac{\nabla f_1^k}{f_1^k}+\frac{1-l^2}{r^2}w_1^k=
\sum_{j=2}^na_{1j}h_j^{k}(0)e^{V^{k}_j}(\frac{w_1^kf_j^k-\psi_j^k}{f_1^k}). $$
Evaluating both sides at $r_k$, the left hand side is strictly positive while the right hand side is non-positive by the definition of $w_1^k$. This contradiction proves (\ref{11jun15e1}). Finally since $g(r)=O(r^2\log\frac 1r)$ near $0$, $g_i^{k,l}(r)=O(\epsilon_k^2r^2)$ near $0$ for $l\ge l_0\ge 3$.
Lemma \ref{11jun14lem3} is established. $\Box$

\medskip

\begin{lem}\label{11jun14lem4} Given $\delta>0$, there exist $C(\delta)>0$ independent of $k,l$ and $k_0(\delta)>1$ such that for $l\ge 3$ and $k\ge k_0$
\begin{equation}\label{11may2e4}
|g^{k,l}_i(r)|\le C\epsilon_k^{(m-2-\delta)^+}(\epsilon_kr)^l+\frac{C}{l^2}\epsilon_k^2r^2(1+r)^{2-m+\delta},\quad r\le \frac 12 \epsilon_k^{-1}
\end{equation}
where
$$(m-2-\delta)^+=\left\{\begin{array}{ll}
m-2-\delta \quad \mbox{ if }\quad m\le 3 \\
1-\delta\quad \mbox{ if }\quad m>3.
\end{array}
\right.
$$
\end{lem}

\noindent{\bf Proof:}

By Lemma \ref{11jun14lem3}, (\ref{11mar31e1}) can be rewritten as
as
$$
(\frac{d^2}{dr^2}+\frac 1r\frac{d}{dr}-\frac{l^2}{r^2})g_i^{k,l}=\bar h_{ikl},\quad
g_i^{k,l}(0)=0
$$
where
\begin{equation}\label{11may2e2}
|\bar h_{ikl}(r)|\le C\epsilon_k^2(1+r)^{2-m+\delta},\quad k\ge k_0(\delta)
\end{equation}
for some $C$ independent of $k$ and $l$.
By standard ODE theory
\begin{equation}\label{11may2e1}
g_i^{k,l}(r)=c^i_{1kl}r^l+c^i_{2kl}r^{-l}-\frac{r^{l}}{2l}\int_r^{\epsilon_k^{-1}}s^{-l+1}\bar h_{ikl}(s)ds
-\frac{r^{-l}}{2l}\int_0^rs^{l+1}\bar h_{ikl}(s)ds.
\end{equation}
When $r\to 0$, it is easy to see that the last two terms in (\ref{11may2e1}) both tend to $0$. Thus $c^i_{2kl}=0$. Let $I_3$ and $I_4$ represent the last two terms in (\ref{11may2e1}), respectively. By (\ref{11may2e2}), we have
\begin{equation}\label{11may2e3}
|I_3(r)|+|I_4(r)|\le \frac{C\epsilon_k^2}{l^2}(1+r)^{4-m+\delta}
\end{equation}
where $C$ is independent of $k$ and $l$. On the other hand, using the information of $w^{1,k}$ for $r\sim \epsilon_k^{-1}$, we know
\begin{equation}\label{12apr27e5}
|g^{k,l}_i(r)|\le C\epsilon_k^{(m-2-\delta)^+}, \quad r\sim \epsilon_k^{-1}.
\end{equation}
In regard to (\ref{11may2e3}) and (\ref{11may2e1}) we have
\begin{equation}\label{12apr27e6}
|c^i_{1kl}|\le C\epsilon_k^{(m-2-\delta)^++l}.
\end{equation}
Hence (\ref{11may2e4}) follows immediately. Lemma \ref{11jun14lem4} is established. $\Box$

\medskip

Let $z^k=(z_1^k,...,z_n^k)$ be the projection of $w^{1,k}$ to
$span\{ \sin l\theta, \cos l\theta, l\ge 2\}$, i.e.
$$z_i^{k}=\sum_{l=2}^{\infty}\bigg (g_i^{k,l}(r)\sin l\theta + \tilde g_i^{k,l}(r)\cos l \theta \bigg )$$
where $\tilde g_i^{k,l}$ is the projection of $w_i^{1,k}$ on $\cos l \theta$. $\tilde g_i^{k,l}$ has similar estimates as that for $g_i^{k,l}$.
Then Lemma \ref{11jun14lem4} leads to
\begin{equation}\label{11jun17e1}
|z_i^k(y)|\le C\epsilon_k^2(1+|y|)^{4-m+\delta}, \quad |y|\le \frac 12\epsilon_k^{-1}
\end{equation}
for $m\le 3$. However for $m>3$, Lemma \ref{11jun14lem4} only gives
\begin{equation}\label{12jun20e5}
|z_i^k(y)|\le C\epsilon_k^{3-\delta}(1+|y|)^2+C\epsilon_k^2(1+|y|)^{4-m+\delta}, \quad |y|\le \frac 12\epsilon_k^{-1}.
\end{equation}
In the following we shall get rid of the first term on the right hand side of (\ref{12jun20e5}). To this end, we need to evaluate the value of $w^{1,k}$ on $\partial \Omega_{1,k}$.

\begin{lem}\label{wikbry}
(a)\,\, If $m<4$, then
\begin{equation}\label{11apr11e10}
|w_i^{1,k}(y)|\le C\epsilon_k^{m^k-2}\log \epsilon_k^{-1},\quad y\in \partial \Omega_{1,k}.
\end{equation}
(b)\,\, If $m=4$ and $|m_i^k-4|\le C/\log \epsilon_k^{-1}$ for all $i$, then
\begin{equation}\label{11jun13e3}
|w_i^{1,k}(y)|\le C\epsilon_k^2(\log \epsilon_k)^2,\quad y\in \partial \Omega_{1,k}.
\end{equation}
\end{lem}

\begin{rem} The assumption $|m_i^k-4|\le C/\log \epsilon_k^{-1}$ when $m=4$ is natural and will be justified in the proof of the main theorems in section 6. We also remark that in (\ref{11apr11e10}) we use $\epsilon_k^{m^k-2}\log 1/\epsilon_k$ instead of the crude $\epsilon_k^{m-2-\delta}$ as before.
\end{rem}

\noindent{\bf Proof:}

The Green's representation formula for $w^{1,k}_i$ gives (see (\ref{newwk}))
\begin{eqnarray}
O(\epsilon_k^2)=w_i^{1,k}(0)=\int_{\Omega_{1,k}}G(0,\eta)\big ( \sum_ja_{ij}h_j^{1,k}(0)e^{V_j^{k}}w_j^{1,k}(\eta)-E_i^k(\eta))d\eta \nonumber\\
-\int_{\partial \Omega_{1,k}}\partial_{\nu}G(0,\eta)w_i^{1,k}(\eta)dS_{\eta} \label{11jun13e5}
\end{eqnarray}
where $G(\cdot, \cdot)$ is the Green's function on $\Omega_{1,k}$ with respect to the Dirichlet boundary condition.
We observe that
\begin{equation}\label{11jun13e6}
G(0,\eta)=-\frac 1{2\pi}\log |\eta |+\frac 1{2\pi} \log \epsilon_k^{-1}+O(\epsilon_k^2).
\end{equation}
Indeed, since $\Omega_{1,k}$ is a translation of $\Omega_k$ by $p_{1,k}$ ( recall $p_{1,k}=O(\epsilon_k)$), the oscillation of $\log |\eta|$ on
$\partial \Omega_{1,k}$ is $O(\epsilon_k^2)$. Thus the oscillation of the regular part of $G(0, \eta)$ is $O(\epsilon_k^2)$, which leads to (\ref{11jun13e6}).
On the other hand let $\bar w_i^{1,k}$ be the average of $w^{1,k}_i$ on $\partial \Omega_{1,k}$, using the fact that $v^{1,k}_i$ is constant on
$\partial \Omega_{1,k}$ and $V_i^k$ has oscillation $O(\epsilon_k^2)$ on $\partial \Omega_{1,k}$ we have
\begin{equation}\label{11jun13e7}
w_i^{1,k}(y)=\bar w_i^{1,k}+O(\epsilon_k^2), \quad \forall y\in \partial \Omega_{1,k}.
\end{equation}
Thus
\begin{equation}\label{11jun13e8}
-\int_{\partial \Omega_{1,k}}\partial_{\nu}G(0,\eta)w_i^{1,k}(\eta)dS_{\eta}
=\bar w_i^{1,k}+O(\epsilon_k^2).
\end{equation}
By using (\ref{11jun13e6}), (\ref{11jun13e7}) and (\ref{11jun13e8}) in (\ref{11jun13e5}) we have
$$
-\bar w_i^{1,k}+O(\epsilon_k^2)=\int_{\Omega_{1,k}}(\frac 1{2\pi}\log \frac{\epsilon_k^{-1}}{|\eta |})\bigg (\sum_j a_{ij}h_j^{k}(0)e^{V_j^{k}}w_j^{1,k}(\eta)
-E_i^k(\eta)\bigg )d\eta
$$
To evaluate the right hand side, we divide $\Omega_{1,k}$ into a symmetric part: $D_1:=B(0, \epsilon_k^{-1}-|p_{1,k}|)$ and a nonsymmetric part:
$\Omega_{1,k}\setminus D_1$ and use $I_1$ and $I_2$ to represent the corresponding integrals on them.
If $\eta\in \Omega_{1,k}\setminus D_1$, it is easy to see from $|p_{1,k}|=O(\epsilon_k)$ that
$$\log (\frac{\epsilon_k^{-1}}{|\eta |})=O(\epsilon_k^2) \mbox{ and } | \Omega_{1,k}\setminus D_1|=O(1). $$
Moreover for $\eta\in \Omega_{1,k}\setminus D_1$, by (\ref{wi0k}) and (\ref{ewik1})
$$w_j^{1,k}(\eta)=O(\epsilon_k^{1-\epsilon}),
\quad E_i^k(\eta)=O(\epsilon_k^{m^k}).$$ Combining these facts we have
$I_2=O(\epsilon_k^2)$.

To evaluate $I_1$, since $\log (\epsilon_k^{-1}/|\eta |)$ is radial, by symmetry only the projections of $w_j^{1,k}$ and $E_i^k$ onto $1$ remain. By (\ref{11mar29e4}) and (\ref{ewik1})
$$|I_1|\le C\log \epsilon_k^{-1}\int_{D_1}\epsilon_k^2(1+|\eta |)^{2-m^k}d\eta. $$
(Note that we use the fact that $e^{V_k(x)}\le C(1+|y|)^{-m^k}$ for $C$ independent of $k$.)
Therefore if $m<4$, $I_1=O(\epsilon_k^{m^k-2})\log \epsilon_k^{-1}$.
If $m=4$ and $|m_i^k-4|\le C\log \epsilon_k^{-1}$ we have  $(1+r)^{-m^k}\le C(1+r)^{-4}$, thus by elementary computation
$$|I_1|\le C\epsilon_k^2(\log \epsilon_k)^2. $$
Lemma \ref{wikbry} is established. $\Box$

\medskip

By Lemma \ref{wikbry}, (\ref{12apr27e5}) can be replaced by
$$
|g^{k,l}_i(r)|\le C\epsilon_k^{m-2-\delta}, \quad r\sim \epsilon_k^{-1}.
$$
Correspondingly, the estimate for $c_{1kl}^i$ becomes
$$|c_{1kl}^i|\le C\epsilon_k^{m-2-\delta+l}. $$
Then it is easy to see that the first term in (\ref{12jun20e5}) can be removed.

\medskip

\noindent{\bf Step five:}

\medskip

For the projection on $span\{\sin \theta, \cos \theta\}$, we
write (\ref{1030e1}) as
$$(\frac{d^2}{dr^2}+\frac 1r\frac{d}{dr}-\frac{1}{r^2})G_{1,i}^k=h(r),\quad 0<r<\epsilon_k^{-1} $$
where $h$ is the collection of other terms. By (\ref{11mar21e1})
$$|h(r)|\le C(1+r)^{1-m+\delta}.$$
Then
\begin{equation}\label{12mar8e1}
G_{1,i}^k=c_{1k}r+\frac{c_{2k}}r-\frac r2\int_r^{\infty}h(s)ds-\frac{r^{-1}}2\int_0^r s^2 h(s)ds.
\end{equation}
Since $G_{1,i}$ is bounded near $0$, $c_{2k}=0$. Using $G_{1,i}^k(\epsilon_k^{-1})=O(\epsilon_k^{m-3-\delta})$ we have
$$|c_{1k}|\le C\epsilon_k^{m-2-\delta}. $$
Then it is easy to see from (\ref{12mar8e1}) that
\begin{equation}\label{11jun23e1}
|G_{t,i}^k(r)|\le Cr(1+r)^{2-m+\delta}, \, t=1,2, \, \, \mbox{ when } m>3.
\end{equation}
Similarly (\ref{11jun23e1}) also holds for $G_{2,i}^k$.

\medskip

Combining the results in the five steps we arrive at the following estimate without distinguishing $m<4$ or not.
\begin{thm}\label{cruexp} Given $\delta>0$, there exist $C(\delta)>0$, $k_0(\delta)>1$ such that for $|y|\le \epsilon_k^{-1}/2$ and $|\alpha |=0,1$, the following holds for all $k\ge k_0$
\begin{equation}\label{11jun14e5}
|D^{\alpha}(v_i^{1,k}(y)-V_i^k(y)-\Phi_i^k(y))|\le C\epsilon_k^2(1+|y|)^{4-m-|\alpha |+\delta}.
\end{equation}
where
$$\Phi_i^k(y)=\epsilon_k(G_{1,i}^k(r)\cos\theta+G_{2,i}^k(r)\sin \theta) $$
with
\begin{equation}\label{11jun14e6}
|G_{t,i}^k(r)|\le Cr(1+r)^{2-m+\delta}\quad t=1,2.
\end{equation}
\end{thm}

Note that the estimate for $|\alpha |=0$ follows directly from the five steps.
The estimate for $|\alpha |=1$ follows from standard gradient estimate for elliptic equations.

\medskip

 Theorem \ref{cruexp} does not distinguish $m<4$ or $m=4$. In the following we apply Theorem \ref{cruexp} to obtain more accurate estimates for $m<4$ and $m=4$, respectively. Both results in the sequel (Theorem \ref{expthm1} and Theorem \ref{expthm2}) play a crucial role in determining the location of maximum points of bubbling solutions.

\begin{thm}\label{expthm1}
Suppose $m<4$, then for
$|y|\le \epsilon_k^{-1}/2$ and $i\in I$
\begin{eqnarray}\label{mle4e1}
&&|D^{\alpha}\big (v_i^k(y)-\phi_i^k(\epsilon_ky)-V_i^k(y-p_{1,k})-\Phi^k_i(y-p_{1,k})\big )|\\
&\le &
C\epsilon_k^2(1+|y|)^{4-m^k-l}\log (2+|y|)
\quad |\alpha |=0,1 \nonumber
\end{eqnarray}
where $v_i^k,\phi_i^k,V_i^k,\Phi_i^k$ are defined by (\ref{vikdef}), (\ref{816e5}), (\ref{88e2}) and (\ref{11aug11e1}), respectively. Moreover
$G_{t,i}^k$ ($t=1,2, i\in I$) satisfy
\begin{equation}\label{11jun9e1}
|G_{t,i}^k(r)|\le Cr(1+r)^{2-m^k}\quad 0<r<\epsilon_k^{-1}.
\end{equation}
\end{thm}

\noindent{\bf Proof:}

We use the same notations as in the proof of Proposition \ref{cruexp}.
First we consider the radial part of $w^{1,k}$: Recall $m<4$. Using Proposition \ref{cruexp} and (\ref{11apr8e1}) we have
\begin{eqnarray*}
&&(\frac{d^2}{dr^2}+\frac 1r\frac{d}{dr})g_i^{k,0}\\
&=&-\sum_j a_{ij}h_j^k(0)e^{V_j^k}g_j^{k,0}-\frac{\epsilon_k^2}4\sum_ja_{ij}\Delta H_j^{1,k}(0)r^2e^{V_j^k}+O(\epsilon_k^2)(1+r)^{\epsilon-m}\\
&=&O(\epsilon_k^2)(1+r)^{2-m^k},\quad 0<r<\epsilon_k^{-1},
\end{eqnarray*}
and
$$ g_i^{k,0}(0)=O(\epsilon_k^2).$$

Multiplying $r$ on both sides of the equation and integrating, we obtain
$$|g_i^{k,0}(r)|\le C\epsilon_k^2\log (2+r). $$

Next we consider the projection of $w^{1,k}$ on $\sin\theta $ and $\cos\theta$.
Let
$\omega^k=(\omega_1^k,...,\omega_n^k)$ be
$$\omega_i^k=w^{1,k}_i-g^{k,0}_i. $$
Then $\omega^k$ satisfies
$$\Delta \omega_i^k+\sum_{j=1}^n a_{ij}h^{k}_j(0)e^{V^{k}_j}\omega_j^k=E_{1,i}^k, \quad \Omega_{1,i} $$
where $E_{1,i}^k$ is the projection of the original right hand side on the subspace spanned by $\sin k\theta$ and $\cos ku$ ($k=1,2,..$).
Since $g_i^{k,0}$ is a radial function, from the asymptotic behavior of $g^{k,0}$ it is easy to see that the oscillation of it on $\partial \Omega_{1,k}$ is $O(\epsilon_k^2)$, therefore the oscillation of $\omega_i^k$ on
$\partial \Omega_{1,k}$ is $O(\epsilon_k^2)$. Let $\bar \omega_i^k$ be the average of $\omega_i^k$ on $\partial \Omega_{1,k}$, then by $\omega_i^k(0)=0$ we have
$$0=\int_{\Omega_{1,k}}G(0,\eta)(\sum_{j=1}^n a_{ij}h^{k}_j(0)e^{V^{k}_j}\omega_j^k(\eta)-E_{1,i}^k(\eta))d\eta+\bar \omega_i^k+O(\epsilon_k^2). $$
Using (\ref{11jun13e6}) in the equation above we have
$\bar \omega_i^k=O(\epsilon_k^2)$, thus
\begin{equation}\label{12mar9e3}
\omega_i^k=O(\epsilon_k^2),\quad \mbox{ on } \quad \partial \Omega_{1,k}.
\end{equation}
Since we have known that $ \epsilon_k G_{t,i}^k(\epsilon_k^{-1})=O(\epsilon_k^{m^k-2})$, we can improve the estimate of $\Phi^k$ using (\ref{12mar8e1}). The estimate of $G_{t,i}^k$ now is
$$|G_{t,i}^k(r)|\le Cr(1+r)^{2-m^k}, \quad r<\epsilon_k^{-1} $$
and
$$|\frac{d}{dr}G_{t,i}^k|\le Cr(1+r)^{1-m^k},\quad 0<r<\epsilon_k^{-1}/2$$
which leads to
\begin{equation}\label{11jun16e1}
|\nabla \Phi^{k}_i(y)|\le C\epsilon_kr(1+r)^{1-m^k},\quad |y|<\epsilon_k^{-1}/2.
\end{equation}
As far as the projection of $w^{1,k}$ on higher frequencies is concerned, since we have (\ref{12mar9e3}), for
$g^{k,l}$ we now have, instead of Lemma \ref{11jun14lem4}
\begin{equation}\label{11jun16e2}
|g^{k,l}_i(r)|\le C\epsilon_k^{m^k-2}(\epsilon_k r)^{l}+\frac{C}{l^2}\epsilon_k^2r^2(1+r)^{2-m^k}, l\ge 2.
\end{equation}
As before we let
$z_i^k=w^{1,k}_i-\Phi_i^{k}$, then $z_i^k$ satisfies
$$\left\{\begin{array}{ll}
\Delta z_i^k=O(\epsilon_k^2)|y|^2(1+|y|)^{-m^k},\quad \Omega_{1,k}\subset \mathbb R^2, \\
\\
z_i^k=O(\epsilon_k^2)\quad \mbox{ on }\quad \partial \Omega_{1,k}.
\end{array}
\right.
$$
Because of (\ref{11jun16e2}) we have
$$|z_i^k(y)|\le C\epsilon_k^2(1+|y|)^{4-m^k}, \quad |y|\le \epsilon_k^{-1}. $$
By standard re-scaling method
\begin{equation}\label{11jun16e3}
|\nabla z_i^k(y)|\le C\epsilon_k^2(1+|y|)^{3-m^k},\quad |y|<\epsilon_k^{-1}/2.
\end{equation}
We have established
$$|D^{\alpha}(v_i^{1,k}(y)-V_i^{k}(y)-\Phi_i^{k}(y)|\le C\epsilon_k^2(1+|y|)^{4-m^k-|\alpha |}\log (2+|y|) $$
for $|\alpha |=0,1$ and $|y|\le \epsilon_k^{-1}/2$.
Recall that $v_i^{1,k}$ is defined in (\ref{vi0kd}). Instead of using the coordinate around $p_{1,k}$ we use the coordinate around the origin to obtain (\ref{mle4e1}).

Theorem \ref{expthm1} is established. $\Box$

\medskip

\begin{thm}\label{expthm2}
If $m=4$ and $|m_i^k-4|\le C/\log \epsilon_k^{-1}$ for all $i\in I$, then we have, for
$|y|\le \epsilon_k^{-1}/2$ and $i\in I$
\begin{eqnarray}\label{m4e1}
&&| D^{\alpha}(v_i^k(y)-\phi_i^k(\epsilon_ky)-V_i^k(y-p_{1,k})-\Phi_i^{k}(y-p_{1,k})\big )|\\
&\le &
C\epsilon_k^2(1+|y|)^{-|\alpha |}(\log (2+|y|))^2.
\quad |\alpha |=0,1, \nonumber
\end{eqnarray}
where $\Phi^{k}$ is of the form stated in (\ref{11aug11e1})
with $G_{t,i}^k$ ($t=1,2$) satisfying
\begin{equation}\label{11jun9e3}
|G_{t,i}^k(r)|\le Cr(1+r)^{-2},\quad 0<r<\epsilon_k^{-1},\quad i\in I.
\end{equation}
\end{thm}

\noindent{\bf Proof:}

For $m=4$ we use (\ref{11mar29e4}) to write (\ref{11apr8e1}) as
$$\left\{\begin{array}{ll}
(\frac{d^2}{dr^2}+\frac 1r\frac{d}{dr})g_i^{k,0}=-\epsilon_k^2\sum_j a_{ij}\frac{\Delta H_j^{1,k}(0)}4 r^2e^{V_j^{k}}+O(\epsilon_k^2)(1+r)^{-4+\epsilon}, \\
\\
g_i^{k,0}(0)=O(\epsilon_k^2),\quad \frac{d}{dr}g_i^{k,0}(0)=0.
\end{array}
\right.
$$
Note that $(1+r)^{-m^k}\le C(1+r)^{-4}$.
Multiplying both sides by $r$ and integrating, we have
$$
g_i^{k,0}(r)=-\epsilon_k^2\sum_j a_{ij}\frac{\Delta H_j^{1,k}(0)}4\int_0^rt^3e^{V_j^{k}(t)}(\log r-\log t)dt+O(\epsilon_k^2)\log (1+r). $$
To evaluate the integral, we use (\ref{12apr24e1}) to write
$$H_i^{1,k}(0)e^{V_i^{k}(t)}=e^{c_i^k}t^{-m^k}+O(t^{-4-\epsilon_0}), \quad t>1 $$
for some $\epsilon_0>0$. Under the assumption of Theorem \ref{expthm2}: $|m_i^k-4|\le C/\log \epsilon_k^{-1}$ we have
$t^{-m^k}=O(t^{-4})$, thus
$$|g_i^{k,0}(r)|\le C\epsilon_k^2(\log (2+r))^2,\quad 0<r<\epsilon_k^{-1}. $$
The projection on higher frequencies has the same estimates as in the case for $m<4$. Specifically, let $\omega_i^k$ be the same as in Theorem \ref{expthm1}. Then (\ref{12mar9e3}) also holds. Correspondingly (\ref{11jun16e1}), (\ref{11jun16e2}) and (\ref{11jun16e3}) still hold with $m^k=4$ by the same proof.
Theorem \ref{expthm2} is established. $\Box$

\section{Location of the blowup points}

In this section we determine the locations of the blowup points in Theorem \ref{expthm1} and Theorem \ref{expthm2}.
\begin{thm}\label{location}
Under the same assumptions as in Theorem \ref{expthm1}
\begin{equation}\label{10050404}
|\sum_i\bigg (\frac{\partial_l h_i^k(0)}{h_i^k(0)}+\partial_l \phi_i^k(0)\bigg )\sigma_i^k|
\le C\epsilon_k^{m^k-2},\quad l=1,2,
\end{equation}
where $C$ is independent of $k$. On the other hand, under the assumptions in Theorem \ref{expthm2} we have
\begin{equation}\label{11jun28e1}
|\sum_i\bigg (\frac{\partial_l h_i^k(0)}{h_i^k(0)}+\partial_l \phi_i^k(0)\bigg )\sigma_i^k|
\le C\epsilon_k^2\log \epsilon_k^{-1},\quad l=1,2,
\end{equation}
where $\sigma_i^k$ is defined in Definition \ref{11def1}.
\end{thm}

\noindent{\bf Proof of Theorem \ref{location}:}
Recall $H_i^k(\epsilon_ky)=h_i^k(\epsilon_ky)e^{\phi_i^k(\epsilon_ky)}$ and $H_i^{1,k}$ is defined in (\ref{hi0k}).
Let $\tilde \Omega_k=B(0, \epsilon_k^{-1}/2)$, we use the following Pohozaev identity for the equation for $v^{1,k}$:
For $\xi\in S^1$,
\begin{eqnarray}\label{10050406}
&&\qquad \sum_i\int_{\tilde \Omega_k}\partial_{\xi}H_i^{1,k}(y)e^{v_i^{1,k}(y)}\\
&=&\int_{\partial \tilde \Omega_k}(\xi\cdot \nu)\sum_ie^{v_i^{1,k}}H_i^{1,k}
+\sum_{ij}a^{ij}\bigg (\partial_{\nu}v_i^{1,k}\partial_{\xi}v_j^{1,k}
-\frac 12\nabla v_i^{1,k}\nabla v_j^{1,k}(\xi \cdot \nu)\bigg )\nonumber
\end{eqnarray}
According to the definition of $H_i^{1,k}$ in (\ref{hi0k})
$$\partial_{\xi}H_i^{1,k}(y)=\epsilon_k\partial_{\xi}H_i^k(0)+\sum_{l=1}^2\epsilon_k^2\partial_{\xi l}H_i^k(0)y^l+O(\epsilon_k^3(1+|y|)^2). $$
Using the expansion of $v_i^{1,k}$ in Proposition \ref{cruexp} ( which holds for $m<4$ and $m=4$) we have
\begin{eqnarray*}
&&\int_{\tilde \Omega_k}\partial_{\xi}H_i^{1,k}(y)e^{v_i^{1,k}(y)}\\
&=&\int_{\tilde \Omega_k}(\epsilon_k\partial_{\xi}H_i^k(0)+\sum_{l=1}^2\epsilon_k^2\partial_{\xi l}H_i^k(0)y^l+O(\epsilon_k^3(1+|y|)^2))\\
&&\cdot \bigg (e^{V_i^{k}(y)}(1+\Phi_i^{k})+O(\epsilon_k^2)(1+|y|)^{4-2m+\epsilon} \bigg )dy.
\end{eqnarray*}
By symmetry we have
$$\int_{\tilde \Omega_k}\sum_{l=1}^2\epsilon_k^2\partial_{\xi l}H_i^k(0)y^l e^{V_i^{k}(y)}
=\int_{\tilde \Omega_k} \epsilon_k\partial_{\xi}H_i^k(0) e^{V_i^{k}(y)} \Phi_i^k(y)dy=0. $$
Also by elementary estimates we have
$$\int_{\tilde \Omega_k}
\big (\sum_{l=1}^2\epsilon_k^2\partial_{\xi l}H_i^k(0)y^l+O(\epsilon_k^3(1+|y|)^2) \big )
\big (e^{V_i^{k}(y)}\Phi_i^{k}+O(\epsilon_k^2)(1+|y|)^{4-2m+\epsilon} \big )=O(\epsilon_k^{m^k-1}) $$
and
$$\int_{\tilde \Omega_k}\epsilon_k\partial_{\xi}H_i^k(0)O(\epsilon_k^2)(1+|y|)^{4-2m+\epsilon}=O(\epsilon_k^{m^k-1}). $$
Thus
\begin{eqnarray*}
&&\int_{\tilde \Omega_k}\partial_{\xi}H_i^{1,k}(y)e^{v_i^{1,k}(y)}\\
&=&\int_{\tilde \Omega_k} \epsilon_k\partial_{\xi}H_i^k(0) e^{V_i^{k}(y)} + O(\epsilon_k^3(1+|y|)^2) e^{V_i^k(y)}.
\end{eqnarray*}

The only difference on whether $m<4$ or $m=4$ is on this term.
For $m<4$, direct computation gives
\begin{eqnarray*}
&&\int_{\tilde \Omega_k}\partial_{\xi}H_i^{1,k}(y)e^{v_i^{1,k}(y)}\\
&=&\epsilon_k\partial_{\xi}H_i^k(0)\int_{\tilde \Omega_k}e^{V_i^{k}(y)}dy
+O(\epsilon_k^{m^k-1}).\\
&=&\epsilon_k\frac{\partial_{\xi}H_i^k(0)}{H_i^k(0)}\int_{\tilde \Omega_k} H_i^k(0)e^{V_i^{k}(y)}dy+O(\epsilon_k^{m^k-1})\\
&=&2\pi \epsilon_k \partial_{\xi}(\log H_i^k)(0)\sigma_i^k+O(\epsilon_k^{m^k-1}).
\end{eqnarray*}
On the other hand for $m=4$, using the closeness between $m_i^k$ and $4$, we have
\begin{eqnarray*}
&&\int_{\tilde \Omega_k}\partial_{\xi}H_i^{1,k}(y)e^{v_i^{1,k}(y)}\\
&=&2\pi \epsilon_k \partial_{\xi} (\log H_i^k)(0)\sigma_i^k+O(\epsilon_k^{3}\log\epsilon_k^{-1}).
\end{eqnarray*}

Next we estimate the right hand side of the Pohozaev identity. The key point in this part is that the only difference between Theorem \ref{expthm1} and Theorem \ref{expthm2} is on the radial part of $w_i^{1,k}$: $g_i^{k,0}$. The projections of $w_i^{1,k}$ on higher frequencies ($ \omega_i^k$ ) have the same estimates. i.e. (\ref{12mar9e3})- (\ref{11jun16e3}) hold for $m<4$ and $m=4$. The detail is as follows:  First by the decay rate of $v_i^{1,k}$ one sees easily that
$$\int_{\partial \tilde \Omega_k}(\xi\cdot \nu)\sum_ie^{v_i^{1,k}}H_i^{1,k}=O(\epsilon_k^{m^k-1}). $$
Before evaluating the remaining two terms we first observe that by (\ref{11jun16e1}) and (\ref{11jun16e3}) that
$$\nabla v_i^{1,k}=\nabla V_i^{k}+\nabla g_i^{k,0}+O(\epsilon_k^{m^k-1}),\quad y\in \partial \tilde \Omega_k. $$
Moreover, by symmetry we have
$$\int_{\partial \tilde \Omega_k}\partial_{\nu}V_i^{k}\partial_{\xi}g^{k,0}_j=\int_{\partial \tilde \Omega_k}\partial_{\xi}V_i^{k}\partial_{\nu}g^{k,0}_j=0. $$
Thus
\begin{eqnarray*}
&&\sum_{i,j}\int_{\partial \tilde \Omega_k}a^{ij}\partial_{\nu}v_i^{1,k}\partial_{\xi}v_j^{1,k} \\
&=&\sum_{i,j}\int_{\partial \tilde \Omega_k}a^{ij}\big ((\partial_{\nu}(V_i^{k}+g_i^{k,0})+O(\epsilon_k^{m^k-1})\big )
\big (\partial_{\xi}(V_j^{k}+g^{k,0}_j)+O(\epsilon_k^{m^k-1})\big )\\
&&=O(\epsilon_k^{m^k-1}).
\end{eqnarray*}

Similarly in the last integral:
\begin{eqnarray*}
&&\sum_{i,j}\int_{\partial \tilde \Omega_k}a^{ij}\nabla v_i^{1,k}\nabla v_j^{1,k} (\xi\cdot \nu)\\
&=&\sum_{i,j}\int_{\partial \tilde \Omega_k}a^{ij}(\frac{d}{dr}(V_i^{k}+g_i^{k,0})+O(\epsilon_k^{m^k-1}))
(\frac{d}{dr}(V_j^{k}+g_j^{k,0})+O(\epsilon_k^{m^k-1}))(\xi\cdot \nu)\\
&=&O(\epsilon_k^{m^k-1}).
\end{eqnarray*}
Theorem \ref{location} is established. $\Box$

\bigskip

\section{The leading term for $\rho^k\to \rho\in \Gamma_1$ }

In this section we complete the proofs of Theorem \ref{rhoto1mle4}, Theorem \ref{rhotome4} and Theorem \ref{locations}. The notation $u^k$ refer to the one in the introduction.
Let
\begin{equation}\label{11jun25e1}
\Theta_i^k=u_i^k-\log \int_M h_ie^{u_i^k}dV_g,\quad i\in I.
\end{equation}
 Then we have
 \begin{equation}\label{1013e1}
 -\Delta_g \Theta_i^k=\sum_{j=1}^na_{ij}\rho_j^k(h_je^{\Theta_j^k}-1).
 \end{equation}
 where
 \begin{equation}\label{1013e2}
 \int_Mh_ie^{\Theta_i^k}dV_g=1,\quad Vol(M)=1.
 \end{equation}
 $\Theta^k=(\Theta_1^k,..,\Theta_n^k)$ is a sequence of blow-up solutions.
 We use $p_k$ to denote the point where the maximum of $\Theta^k$ on $M$ is taken. Take the local coordinates around $p_k$. Then $ds^2$ has the form $e^{\psi(y_{p_k})}(dy_1^2+dy_2^2)$ where
 $$|\nabla\psi(0)|=0,\quad \psi(0)=0,\quad \Delta \psi=-2Ke^{\psi},$$
  $K$ is the Gauss curvature. In local coordinates, (\ref{1013e1}) becomes
 \begin{equation}\label{1013e3}
 -\Delta \Theta_i^k=\sum_{j=1}^na_{ij}\rho_j^ke^{\psi}(h_je^{\Theta_j^k}-1), \quad \mbox{in}\quad B_{\delta_0}
 \end{equation}
for some $\delta_0$ small.
 Let $f_i^k$ solve
 \begin{equation}\label{1013e4}\left\{\begin{array}{ll}
 \Delta f_i^k=\sum_{j=1}^na_{ij}\rho_j^ke^{\psi},\quad B_{\delta_0}\\
 \\
 f_i^k(0)=0,\quad \nabla f_i^k(0)=0.
 \end{array}
 \right.
 \end{equation}
and $f_i^k$ has bounded oscillation on $\partial B_{\delta_0}$.
Set
$$\tilde \Theta_i^k=\Theta_i^k(p_k+\cdot)-f_i^k\quad \mbox{and }\quad \tilde h_i^k=\rho_i^k h_i(p_k+\cdot)e^{\psi+f_i^k}, $$
 then $\tilde \Theta_i^k$ satisfies
\begin{equation}\label{1013e5}
-\Delta \tilde \Theta_i^k=\sum_ja_{ij}\tilde h_j^ke^{\tilde \Theta_j^k},\quad B_{\delta_0}.
\end{equation}
Note that from the definition of $\tilde h_i^k$, $f_i^k$ and $\psi$, the following can be verified by direct computation:
\begin{equation}\label{11jun18e1}
\left\{\begin{array}{ll}
\tilde h_i^k(0)=\rho_i^k h_i(p_k), \quad \nabla \tilde h_i^k(0)=\rho_i^k \nabla h_i(p_k), \\
\\
\Delta \tilde h_i^k(0)=\rho_i^k \bigg (\Delta h_i(p_k)+h_i(p_k)(-2K(p_k)+\sum_{j=1}^n a_{ij}\rho_j^k)\bigg ).
\end{array}
\right.
\end{equation}

Let $M_k=\max_{i\in I}\max_M\Theta_i^k$, $\epsilon_k=e^{-\frac 12 M_k}$,
\begin{equation}\label{11mar7e2}
\tilde v_i^k(y)=\tilde \Theta_i^k(\epsilon_ky)+2\log \epsilon_k=\tilde \Theta_i^k(\epsilon_ky)-M_k,
\end{equation}
then we have
\begin{equation}\label{11jun11e1}
\Delta \tilde v_i^k(y)+\sum_{j=1}^n a_{ij} \tilde h_j^k(\epsilon_ky)e^{\tilde v_j^k(y)}=0,\quad |y|\le \delta_0\epsilon_k^{-1}.
\end{equation}

The following lemma proves that the $\Theta^k$ is a fully blown up sequence.
\begin{lem}\label{ful}
Along a subsequence $(\tilde v_1^k,...,\tilde v_n^k)$ converge in $C^2_{loc}(\mathbb R^2)$ to a system of $n$ equations.
\end{lem}

\noindent{\bf Proof:}
By way of contradiction we assume that only $l$ ($l<n$) components of $(\tilde v_1^k,..,\tilde v_n^k)$ converge to a system of $l$ equations. Without loss of generality we assume that the first $l$ components of $(\tilde v_1^k,..,\tilde v_n^k)$ converge to $(\tilde v_1,..,\tilde v_l)$ that satisfies
$$-\Delta \tilde v_i=\sum_{j=1}^l a_{ij} \tilde h_j e^{\tilde v_j}, \quad i=1,..,l. \mbox{ in } \mathbb R^2 $$
where $\tilde h_i=\rho_i h_i(p)$. Let
$$\sigma_{i,v}=\frac 1{2\pi} \int_{\mathbb R^2} \tilde h_i e^{\tilde v_i}. $$
Then the entire solution $(\tilde v_1,...,\tilde v_l)$ with finite energy satisfies
$$\sum_{j=1}^l a_{ij}\sigma_{j, v}>2. $$
Let $J=\{1,...,l\}$, then by Theorem C in \cite{linzhang2} $(\tilde v_1,...,\tilde v_l)$ also satisfies
$$\frac 1{4\pi^2} \Lambda_J(\sigma_v):=4\sum_{i\in J}\sigma_{i,v}-\sum_{i,j\in J} a_{ij} \sigma_{i,v}\sigma_{j, v}=0. $$
By (\ref{1013e2})
$$\frac{\rho_i^k}{2\pi}\ge \frac 1{2\pi}\int_{B(p_k,\delta)}\rho_i^k h_i e^{\Theta_i^k}dV_g
=\frac{1}{2\pi }\int_{B_{\delta}}\tilde h_i^k e^{\tilde \Theta_i^k}dx. $$
Thus by letting $k\to \infty$ we have
$$\frac{\rho_i}{2\pi }\ge \sigma_{i,v},\quad i=1,..,l. $$
Let $s_i=\frac{\rho_i}{2\pi }-\sigma_{i,v}$, $i\in J$, then easy to see
\begin{equation}\label{12jun22e1}
\frac{1}{4\pi^2}(\Lambda_J(\rho)-\Lambda_J(\sigma_v))=-2\sum_{i\in J}(\sum_{j\in J} a_{ij}\sigma_{j,v}-2)s_i-\sum_{i,j}a_{ij}s_is_j\le 0.
\end{equation}
Since $\Lambda_J(\sigma_v)=0$, (\ref{12jun22e1}) is a violation of the definition of $\Gamma_1$ in the introduction. Lemma \ref{ful} is established. $\Box$

\medskip

Let $\phi^k=(\phi_1^k,...,\phi_n^k)$ be the harmonic function that takes $0$ at $0$ and that kills the oscillation of $\tilde \Theta^k$ on
$\partial B_{\delta_0}$.
The first term in the expansion of
$\tilde v^k-\phi^k(\epsilon_k\cdot)$ is $U^k(\cdot -p_{1,k})$ where $U^k=(U_1^k,..,U_n^k)$ satisfies
\begin{equation}\label{11jun25e2}
\left\{\begin{array}{ll}
-\Delta U_i^k=\sum_ja_{ij}\tilde h_j^k(0)e^{U_j^k}\quad \mathbb R^2, \\
U_i^k(0)=\tilde v_i^k(0),\quad i\in I.
\end{array}
\right.
\end{equation}
Moreover, $U_i^k+\log (\rho_i^k h_i(p_k))$ satisfies
$$-\Delta (U_i^k +\log (\rho_i^k h_i(p_k)))=\sum_ja_{ij}e^{U_j^k+\log (\rho_j^k h_j(p_k))}\quad \mathbb R^2. $$
Set
\begin{equation}\label{12apr29e2}
\bar \sigma_i^k=\frac 1{2\pi}\int_{\mathbb R^2}\rho_i^k h_i(p_k) e^{U_i^k},\quad m_i^k=\sum_{j=1}^n a_{ij}\bar \sigma_j^k.
\end{equation}
Correspondingly
$$m^k=\min\{m_1^k,..,m_n^k\},\quad m=\lim_{k\to \infty} m^k. $$
Later we shall show that $m$ is equal to the one defined in (\ref{12apr23e1}).
In regard to Lemma \ref{globalcpt} we set $c_i^k$ as
\begin{eqnarray}\label{11jun25e3}
c_i^k=U_i^k(0)+\log (\rho_i^k h_i(p_k))+\frac 1{2\pi}\int_{\mathbb R^2}\log |\eta |\sum_{j=1}^na_{ij}\tilde h_j^k(0)e^{U_j^k}d\eta
\\
=\Theta_i^k(p_k)-M_k+\log (\rho_i^k h_i(p_k))+\frac 1{2\pi}\int_{\mathbb R^2}\log |\eta |\sum_{j=1}^na_{ij}\rho_j^k h_j(p_k)e^{U_j^k}d\eta,
\nonumber
\end{eqnarray}
and
\begin{equation}\label{12jun29e1}
c_i=\lim_{k\to \infty} c_i^k.
\end{equation}
According to Proposition \ref{cruexp}, the expansion of $\tilde v_i^k$ can be written as
\begin{equation}\label{11jun17e8}
\tilde v_i^k(y)=U_i^k(y-p_{1,k})+\phi_i^k(\epsilon_ky)+\Phi_i^{k}(y-p_{1,k})+O(\epsilon_k^2)(1+|y|)^{4-m+\epsilon}
\end{equation}
on $B(0,\delta_0\epsilon_k^{-1})$, where $p_{1,k}=O(\epsilon_k)$,
$$\Phi_i^{k}(y)=\epsilon_k(G_{1,i}^k(r)\cos \theta +G_{2,i}^k(r)\sin \theta )$$
such that
$$|G_{t,i}^k(r)|\le Cr(1+r)^{2-m^k},\quad 0<r<\delta_0 \epsilon_k^{-1},\quad t=1,2. $$

 Here we note that the Green's function is of the form
$$G(x,\eta)=-\frac 1{2\pi}\log |x-\eta |\chi +\gamma(x,\eta)$$
where $\chi$ is a cut-off function such that $\chi\equiv 1$ in $B_{\epsilon_1}$ for some $\epsilon_1>0$ ($2\epsilon_1$ is less than the
 injectivity radius of $M$) and $\chi\equiv 0$ outside $B_{2\epsilon_1}$.
In the sequel we always assume $\delta_0<\epsilon_1/4$.

The Green's representation formula for $\Theta_i^k$ is
\begin{equation}\label{11may23e1}
\Theta_i^k(x)=\bar \Theta_i^k+\int_MG(x,\eta)(\sum_ja_{ij}\rho_j^kh_je^{\Theta_j^k})dV_g.
\end{equation}

Next we claim that
\begin{eqnarray}\label{nov2e5}
&&|\bar \Theta_i^k+\frac{m_i^k-2}2M_k-c_i^k+\log (\rho_i^k h_i(p_k))+2\pi\gamma(p_k,p_k)m_i^k|\nonumber \\
&\le & C_0(\delta_0)\epsilon_k^{m-2-\delta}
\end{eqnarray}
where $\delta>0$ is small positive number, $\bar \Theta_i^k$ is the average of $\Theta_i^k$ on $M$, $C_0(\cdot)$ is a positive function such that $C_0(\delta_0)\to \infty$ as $\delta_0\to 0$.
To see this, let $x=p_k$ in (\ref{11may23e1}),
\begin{eqnarray}
\Theta_i^k(p_k)&=&\bar \Theta_i^k-\frac 1{2\pi}\int_{B(p_k,\delta_0)}\log |p_k-\eta |\sum_{j}a_{ij}\rho_j^k h_j e^{\Theta_j^k}dV_g \nonumber\\
&&+\int_{B(p_k,\delta_0)}\gamma(p_k, \eta)\sum_{j}a_{ij}\rho_j^k h_j e^{\Theta_j^k}dV_g \nonumber\\
\label{11mar7e0}
&&+\int_{M\setminus B(p_k,\delta_0)}G(p_k,\eta)\sum_ja_{ij}\rho_j^kh_je^{\Theta_j^k}dV_g.
\end{eqnarray}
Observe that
\begin{equation}\label{11mar7e1}
\rho_j^kh_je^{\Theta_j^k}dV_g(\eta)=\rho_j^kh_je^{f_j^k+\psi}e^{\Theta_j^k-f_j^k}d\eta=\tilde h_j^ke^{\tilde \Theta_j^k}d\eta.
\end{equation}
Then by the definition of $\epsilon_k$, $\tilde v^k$, $c_i^k$ we obtain
\begin{eqnarray}\label{11jun24e2}
&&-\frac 1{2\pi}\int_{B(p_k,\delta_0)}\log |p_k-\eta |\sum_ja_{ij}\rho_j^kh_je^{\Theta_j^k}dV_g\\
&=&-\frac 1{2\pi}\int_{B(0,\delta_0\epsilon_k^{-1})}(\log (\epsilon_k|\eta_1|))\sum_ja_{ij}(\tilde h_j^ke^{\phi_j^k})(\epsilon_k\eta_1)e^{\tilde v_j^k(\eta_1)-\phi_j^k(\epsilon_k\eta_1)}d\eta_1 \nonumber\\
&=&-\frac 1{2\pi}\int_{B(0,\delta_0\epsilon_k^{-1})}(\log \epsilon_k+\log |\eta_1|)\sum_ja_{ij}H_j^{1,k}(\eta)
e^{v_j^{1,k}(\eta_1)}
d\eta_1 \nonumber
\end{eqnarray}
where, by the same kind of notations used in section \ref{blowup} and section 4,
$$v_i^{1,k}(\cdot)=\tilde v_i^k(p_{1,k}+\cdot)-\phi_i^k(\epsilon_kp_{1,k}+\epsilon_k\cdot) $$
and
\begin{equation}\label{11jun25e4}
H_i^{1,k}(\cdot )=\tilde h_i^k(\epsilon_k p_{1,k}+\epsilon_k \cdot )e^{\phi_i^k(\epsilon_kp_{1,k}+\epsilon_k \cdot)}.
\end{equation}

By (\ref{11jun18e1}) we have
$$H_i^{1,k}(0)=\rho_i^k h_i(p_k)+O(\epsilon_k^2), \quad \nabla H_i^{1,k}(0)=\epsilon_k \rho_i^k \nabla h_i(p_k)+O(\epsilon_k^3)$$
and
$$H_i^{1,k}(\eta)=\rho_i^k h_i(p_k)+ \epsilon_k \rho_i^k \nabla h_i(p_k)\cdot \eta +O(\epsilon_k^2)(1+|\eta |)^2. $$
On the other hand by Proposition \ref{cruexp}
$$e^{v_i^{1,k}(\eta_1)}=e^{U_i^{k}}(1+\Phi_i^{k}(\eta_1)+O(\epsilon_k^2)(1+|\eta_1|)^{4-m+\delta}). $$
Hence
\begin{eqnarray}\label{11jun24e1}
&&H_i^{1,k}(\eta)e^{v_i^{1,k}(\eta)}\\
&=& (\rho_i^k h_i(p_k)+\epsilon_k \rho_i^k \nabla h_i(p_k)\cdot \eta +\Phi_i^{k}(\eta))e^{U_i^k(\eta)}
+O(\epsilon_k^2)(1+|\eta |)^{2-m+\delta}. \nonumber
\end{eqnarray}
Using (\ref{11jun24e1}) in the evaluation of (\ref{11jun24e2}) we have
\begin{eqnarray*}
&&-\frac 1{2\pi}\int_{B(p_k,\delta_0)}\log |p_k-\eta |\sum_ja_{ij}\rho_j^kh_je^{\Theta_j^k}dV_g\\
&=&\frac{M_k}2m_i^k+\Theta_i^k(p_k)-M_k-c_i^k+\log (\rho_i^k h_i(p_k))+C_0(\delta_0)(\epsilon_k^{m-2-\delta}).
\end{eqnarray*}

Similarly
$$\int_{B(p_k,\delta_0)}\gamma(p_k,\eta)\sum_ja_{ij}\rho_j^kh_je^{\Theta_j^k}dV_g=\gamma(p_k,p_k)2\pi m_i^k+C_0(\delta_0)(\epsilon_k^{m-2-\delta}). $$
For the final term in (\ref{11mar7e0}) by the following crude estimate established in \cite{linzhang2}
$$\Theta_i^k(x)=-\frac{m_i^k-2}2 M_k +O(1), \quad x\in M\setminus B(p_k, \delta_0), \quad O(1)\to \infty \mbox{ if } \delta_0\to 0 $$
it is easy to see
$$\int_{M\setminus B(p_k,\delta_0)}G(p_k,\eta)\sum_j a_{ij} \rho_j^k h_j e^{\Theta_j^k} dV_g=C_0(\delta_0)O(\epsilon_k^{m-2-\delta}).$$
Thus, back to (\ref{11mar7e0}),
\begin{eqnarray*}
\Theta_i^k(p_k)=\bar \Theta_i^k+\frac{M_k}2m_i^k+\Theta_i^k(p_k)-M_k-c_i^k+\log (\rho_i^k h_i(p_k))\\
+\gamma(p_k,p_k)2\pi m_i^k+O(\epsilon_k^{m-2-\delta}).
\end{eqnarray*}
Then (\ref{nov2e5}) follows.

\medskip

Next for $x\in M\setminus B(p_k,\delta_0)$, by (\ref{11may23e1}) and standard estimates
\begin{equation}\label{11july6e3}
\Theta_i^k(x)=\bar \Theta_i^k+2\pi G(x,p_k)m_i^k+C_0(\delta_0)\epsilon_k^{m-2-\delta}.
\end{equation}
Consequently by (\ref{11july6e3}) and (\ref{nov2e5})
\begin{eqnarray}\label{1216e1}
e^{\Theta_i^k(x)}&=&e^{\bar \Theta_i^k+2\pi G(x,p_k)m_i^k}+E_{\delta_0},\quad x\in M\setminus B(p_k,\delta_0) \\
&=&\epsilon_k^{m_i^k-2}\frac{e^{c_i^k}}{\rho_i^k h_i(p_k)} e^{2\pi m_i^k (G(x,p_k)-\gamma(p_k,p_k))} +E_{\delta_0}. \nonumber
\end{eqnarray}
where
\begin{equation}\label{11131e1}
|E_{\delta_0}|\le C_0(\delta_0)\epsilon_k^{m^k-2+\delta}.
\end{equation}
In the sequel we shall always use $E_{\delta}$ to represent a term bounded by the right hand side of (\ref{11131e1}).

\subsection{Proof of Theorem \ref{rhoto1mle4}.}

  In this case $m<4$.
 Since $\int_M h_ie^{\Theta_i^k}dV_g=1$ we write
$$\rho_i^k=\int_{B(p_k,\delta_0)}\rho_i^kh_ie^{\Theta_i^k}dV_g+\int_{M\setminus B(p_k,\delta_0)}\rho_i^kh_ie^{\Theta_i^k}dV_g
=\rho_{ia}^k+\rho_{ib}^k. $$
By (\ref{11mar7e1}), (\ref{11mar7e2})
$$
\rho_{ia}^k=\int_{B(0,\delta_0)}\tilde h_i^ke^{\tilde \Theta_i^k}d\eta=
\int_{B(0,\delta_0\epsilon_k^{-1})}\tilde h_i^k(\epsilon_ky)e^{\tilde v_i^k(y)}dy.
$$
Let
$$I_1=\{i\in I;\quad \lim_{k\to \infty}m_i^k=\lim_{k\to \infty} m^k.\quad \}. $$
Using the expansion of $\tilde v_i^k$ in (\ref{11jun17e8})  we have (since $m<4$)
\begin{equation}\label{rhoia}
\rho_{ia}^k=\int_{B(0,\delta_0\epsilon_k^{-1})}\tilde h_i^k(0)e^{U_i^k(y)}dy+o(\delta_0)\epsilon_k^{m^k-2},\quad i\in I_1.
\end{equation}
Now for $i\not \in I_1$ we have
\begin{equation}\label{rhoia1}
\rho_{ia}^k=\int_{B(0,\delta_0\epsilon_k^{-1})}\tilde h_i^k(0)e^{U_i^k(y)}dy+E_{\delta_0},\quad i\not \in I_1
\end{equation}
and
\begin{equation}\label{rhobc}
|\rho_{ib}^k|=E_{\delta_0},\quad i\not \in I_1.
\end{equation}
It is easy to see from (\ref{12apr29e2}) and (\ref{11jun18e1}) that the $m$ defined by (\ref{12apr29e2}) is the same as the one in (\ref{12apr23e1}).

Combining (\ref{rhoia}), (\ref{rhoia1}), (\ref{rhobc}) and  Lemma \ref{globalcpt} we have
\begin{eqnarray*}
&&\sum_i \frac{4}{2\pi}\rho_{ia}^k-\sum_{ij}a_{ij}\frac{\rho_{ia}^k}{2\pi}\frac{\rho_{ja}^k}{2\pi}\\
&=&2\sum_{i\in I_1}e^{c_i^k}\delta_0^{2-m_i^k}\epsilon_k^{m_i^k-2}+o(\delta_0)\epsilon_k^{m^k-2}+E_{\delta_0}.
\end{eqnarray*}
Using (\ref{rhoia1}), (\ref{rhobc}) to change from $\rho_{ia}^k$ to $\rho_i^k$,
\begin{eqnarray}\label{1105e2}
&&\frac{4}{2\pi}\sum_i\rho_i^k =\sum_{ij}a_{ij}\frac{\rho_i^k}{2\pi}\frac{\rho_j^k}{2\pi}+2\sum_{i\in I_1}e^{c_i^k}\delta_0^{2-m_i^k}\epsilon_k^{m_i^k-2}\\
&&-2
\sum_{i\in I_1}(m_i^k-2)\frac{\rho_{ib}^k}{2\pi}+o(\delta_0)\epsilon_k^{m^k-2}+E_{\delta_0}. \nonumber
\end{eqnarray}

Using (\ref{1216e1}) we have

\begin{eqnarray*}
\rho_{ib}^k&=&\int_{M\setminus B(p_k,\delta_0)}\rho_i^kh_ie^{\Theta_i^k}dV_g\\
&=&\epsilon_k^{m_i^k-2}\int_{M\setminus B(p_k,\delta_0)}\frac{h_i(x)e^{c_i^k}}{h_i(p_k)}e^{2\pi m_i^k(G(x,p_k)-\gamma(p_k,p_k))}dV_g
+E_{\delta_0},\quad i\in I_1.
\end{eqnarray*}

Combining terms we have
\begin{eqnarray}
&&\frac{4}{2\pi}\sum_i\rho_i^k-\sum_{i,j}a_{ij}\frac{\rho_i^k}{2\pi}\frac{\rho_j^k}{2\pi} \nonumber\\
&=&2\sum_{i\in I_1}e^{c_i^k}\epsilon_k^{m_i^k-2}\bigg (\delta_0^{2-m_i^k}-\frac {m_i^k-2}{2\pi}\int_{M\setminus B(p_k,\delta_0)}
\frac{h_i(x)}{h_i(p_k)}e^{2\pi m_i^k(G(x,p_k)-\gamma(p_k,p_k))}dV_g \bigg )\nonumber\\
&&\qquad \qquad \qquad +o(\delta_0)\epsilon_k^{m^k-2}+E_{\delta_0}. \label{11jun22e1}
\end{eqnarray}

We claim that for fixed $k$ the following limit exists:
\begin{equation}\label{leadterm1}
\lim_{\delta_0\to 0}\bigg (\delta_0^{2-m_i^k}-\frac {m_i^k-2}{2\pi}\int_{M\setminus B(p_k,\delta_0)}
e^{2\pi m_i^k(G(x,p_k)-\gamma(p_k,p_k))}\frac{h_i}{h_i(p_k)}dV_g \bigg )
\end{equation}
Indeed, write the second integral as the sum of one integral over $B(p_k,\delta_1)\setminus B(p_k,\delta_0)$ and the other over $M\setminus B(p_k,\delta_1)$, where $\delta_1$ is chosen small enough so that in $B(p_k,\delta_1)$,
\begin{equation}\label{11mar9e5}
dV_g=e^{\psi}dx\quad \mbox{ where }\psi(x)=O(|x|^2)\mbox{ in }B_{\delta_1}.
\end{equation}
In the local coordinate at $p_k$, it suffices to prove
$$
\lim_{\delta_0\to 0}\bigg (\delta_0^{2-m_i^k}-\frac {m_i^k-2}{2\pi}\int_{B_{\delta_1}\setminus B_{\delta_0}}
e^{2\pi m_i^k(G(p_k+x,p_k)-\gamma(p_k,p_k))}\frac{h_i(p_k+x)}{h_i(p_k)}dV_g(x) \bigg )
$$
exists.
Since
$$G(p_k+x,p_k)=-\frac 1{2\pi}\log |x|+\gamma(p_k+x,p_k),$$
we use (\ref{11mar9e5}) and the Taylor expansions of $h_i$ and $\gamma(p_k+x,p_k)-\gamma(p_k,p_k)$ to obtain
\begin{eqnarray*}
&&e^{2\pi m_i^k(G(p_k+x,p_k)-\gamma(p_k,p_k))}h_i(x+p_k)/h_i(p_k) dV_g\\
&=&|x|^{-m_i^k}(1+\sum_{l=1}^2c_lx^l+O(|x|^2)) dx.
\end{eqnarray*}
where $c_1,c_2$ are some constants. Observe that the terms with $c_1,c_2$ disappear in the integration and $m_i^k<4-\epsilon_1$ for some $\epsilon_1>0$ for all $i\in I_1$, thus the limits in (\ref{leadterm1}) exists.

\medskip

The following lemma states the closeness between $\rho^k$ to $\rho$ if $m<4$. It will be used to simplify the leading terms in the statements of main theorems.

\begin{lem}\label{close1}
Let $\rho^k\to \rho\in \Gamma_1$ ($\rho\neq Q$) such that all $\rho_i^k-\rho_i$ have the same sign,
then
\begin{equation}\label{11jun22e2}
\rho_i^k-\rho_i=O(\epsilon_k^{m-2-\epsilon})\quad \mbox{ and }\quad m_i^k-m_i=O(\epsilon_k^{m-2-\epsilon}).
\end{equation}
\end{lem}

\noindent{\bf Proof of Lemma \ref{close1}:}
Let $s_i^k=\rho_i-\rho_i^k$. Then all $s_i^k$ are of the same sign. By (\ref{11jun22e1}) we have
$$\sum_i (m_i-2)s_i^k=O(\epsilon_k^{m-2-\epsilon}). $$
Thus (\ref{11jun22e2}) holds and Lemma \ref{close1} is established. $\Box$

\medskip

Using (\ref{11jun22e2}) in (\ref{11jun22e1}) we can rewrite the leading term as
\begin{eqnarray*}
&&\frac{4}{2\pi}\sum_i\rho_i^k-\sum_{i,j}a_{ij}\frac{\rho_i^k}{2\pi}\frac{\rho_j^k}{2\pi} \nonumber\\
&=&2\sum_{i\in I_1}e^{c_i^k}\epsilon_k^{m-2}\bigg (\delta_0^{2-m}-\frac {m-2}{2\pi}\int_{M\setminus B(p_k,\delta_0)}
\frac{h_i(x)}{h_i(p_k)}e^{2\pi m(G(x,p_k)-\gamma(p_k,p_k))}dV_g \bigg )\nonumber\\
&&\qquad +o(\delta_0)\epsilon_k^{m-2}+E_{\delta_0}.
\end{eqnarray*}
Since $c_i^k\to c_i$, the $c_i^k$ can be replaced by $c_i$ in (\ref{12jun29e1}).
Theorem \ref{rhoto1mle4} is established. $\Box$

\medskip

Next we establish the following lemma regardless of $m=4$ or not.
\begin{lem}\label{graest}
\begin{equation}\label{graphik}
\nabla \phi_i^k(0)=2\pi m_i^k\nabla_1 \gamma(p_k,p_k)+O(\epsilon_k^{m-2-\epsilon}).
\end{equation}
where $\nabla_1$ means the differentiation with respect to the first component.
\end{lem}

\noindent{\bf Proof of Lemma \ref{graest}:}

First from (\ref{11july6e3}),
for $x\in \partial B(p_k,\delta_0)$
\begin{eqnarray*}
&& \tilde \Theta_i^k(x)=\Theta_i^k(x)-f_i^k(x)\\
&=&\bar \Theta_i^k-m_i^k\log |x-p_k|+2\pi m_i^k \gamma(x,p_k)-f_i^k(x)+O(\epsilon_k^{m-2-\epsilon}).
\end{eqnarray*}
For the first derivative, by standard estimates
$$|D^l(\tilde \Theta_i^k)(x)-D^l(-m_i^k\log |x-p_k|+2\pi m_i^k \gamma(x,p_k)-f_i^k(x))|\le C\epsilon_k^{m-2-\epsilon}, \quad l=0,1$$
for $x\in B(p_k, \delta_0)$.  On the  other hand, we
recall that $\tilde v_i^k(y)=\tilde \Theta_i^k(\epsilon_ky)+2\log \epsilon_k$. According to Proposition \ref{cruexp}
$$|D^l(\tilde v_i^k(y)-U_i^k(y-p_{1,k})-\phi_i^k(\epsilon_k y)-\Phi_i^{k}(y-p_{1,k})|\le C\epsilon_k^2(1+|y|)^{4-m+\epsilon-l}$$
for $l=0,1$ and $|y|\le 2\delta_0\epsilon_k^{-1}$.
Using the asymptotic behavior $U_i^k$ for $|y|\sim \epsilon_k^{-1}$ (\ref{12apr29e1}) we have
$$U_i^k(y-p_{1,k})+\log (\rho_i^kh_i(p_k))=-\frac{m_i^k}2M_k-m_i^k\log |x|+c_i^k+O(\epsilon_k^{m^k-2}) $$
where $|x|=\epsilon_k |y|$. Thus by (\ref{nov2e5}) for $l=0$ we have
\begin{equation}\label{11july7e11}
\phi_i^k(x)=2\pi m_i^k(\gamma(x,p_k)-\gamma(p_k,p_k))-f_i^k(x)+O(\epsilon_k^{2-\epsilon}), \quad x\in B(p_k,2\delta_0)\setminus B(p_k, \delta_0/2).
\end{equation}
Thus the comparison on $l=1$ and (\ref{1013e4}) yield (\ref{graphik}).
Lemma \ref{graest} is established. $\Box$

\medskip

\subsection{Proof of Theorem \ref{rhotome4}}

In this case $m=4$, we first give a rough estimate of $\rho_i^k-\rho_i$.

\begin{lem}\label{lemq4} Let $\rho^k$ tend to $Q$ such that all $\rho_i^k-Q_i$ have the same sign. Then
\begin{equation}\label{11jun17e9}
m_i^k-4=O(\epsilon_k^{2-\epsilon}),\quad \rho_i^k-\rho_i=O(\epsilon_k^{2-\epsilon})\quad \forall i\in I.
\end{equation}
\end{lem}

\noindent{\bf Proof of Lemma \ref{lemq4}:}
Since $\rho^k\to Q$, $m=4$, then
all $m_i^k\to 4$. Recall that $\rho_i^k=2\pi\bar \sigma_i^k+O(\epsilon_k^{2-\epsilon})$,
The rest of the proof is the same as that in the proof of Lemma \ref{close1}. Lemma \ref{lemq4} is established. $\Box$

\medskip

By (\ref{11jun28e1})
$$\sum_i(\partial_l \log \tilde h_i^k(0)+\partial_l \phi_i^k(p_k))\bar \sigma_i^k=O(\epsilon_k^2\log \epsilon_k^{-1}), $$
then by (\ref{11jun18e1})
\begin{equation}\label{locationm4}
\sum_{i=1}^n \bigg (\partial_l(\log h_i+\phi_i^k)(p_k) \bigg )\bar \sigma_i^k=O(\epsilon_k^2)\log \epsilon_k^{-1},\quad l=1,2.
\end{equation}

The distance from $\rho^k$ to $\Gamma_1$ can be computed as follows:
$$\rho_i^k=\int_{B(p_k,\delta_0)}\rho_i^k h_i e^{\Theta_i^k} dV_g +\int_{M\setminus B(p_k,\delta_0)} \rho_i^k h_i e^{\Theta_i^k} dV_g. $$
By (\ref{11jun17e9}) and Theorem \ref{expthm2} the second integral is $O(\epsilon_k^2)$, this is the same as the computation for the single equation \cite{ChenLin1}. Therefore
\begin{eqnarray}\label{11july5e1}
&&\rho_i^k= \int_{B(p_k,\delta_0)}\rho_i^k h_i e^{\Theta_i^k} dV_g+O(\epsilon_k^2)\\
&=&\int_{B(0,\delta_0)}\tilde h_i^k e^{\phi_i^k}e^{\tilde \Theta_i^k-\phi_i^k}d\eta+O(\epsilon_k^2)\nonumber\\
&=&\int_{B(0,\delta_0 \epsilon_k^{-1})}H_i^{1,k}(\eta)e^{v_i^{1,k}(\eta)}d\eta+O(\epsilon_k^2)\nonumber\\
&=&\int_{B(0,\delta_0 \epsilon_k^{-1})} \rho_i^k h_i(p_k) e^{U_i^k(\eta)}d\eta+\int_{B(0,\delta_0 \epsilon_k^{-1})}
\frac 14 \Delta H_i^{1,k}(0)|\eta |^2 e^{U_i^k(\eta)}d\eta \nonumber\\
&&\qquad +O(\epsilon_k^2).\nonumber
\end{eqnarray}
The first integral on the right hand side of the above is $2\pi \bar \sigma_i^k+O(\epsilon_k^2)$. To evaluate the second term on the right hand side, we first use the
definition of the $H_i^{1,k}$ in (\ref{11jun25e4}), (\ref{11jun17e9}), (\ref{11jun18e1}) and (\ref{graphik}) to have
\begin{eqnarray}\label{h1kdel}
\Delta H_i^{1,k}(0)=\rho_i^k h_i(p_k)\bigg (\frac{\Delta h_i(p_k)}{h_i(p_k)}-2K(p_k)+8\pi\\
+16\pi \frac{\nabla h_i(p_k)}{h_i(p_k)}\cdot \nabla_1\gamma(p_k,p_k)
+64\pi^2 |\nabla_1\gamma(p_k,p_k)|^2\bigg )\epsilon_k^2+O(\epsilon_k^{4-\epsilon}). \nonumber
\end{eqnarray}
For $e^{U_i^k}$ we use the definition of $c_i^k$ in (\ref{11jun25e3}) and (\ref{12apr29e1}) to have
\begin{equation}\label{11jun25e8}
\rho_i^k h_i(p_k) e^{U_i^k(\eta)}=e^{c_i^k} |\eta |^{-4}+O(|\eta |^{-4-\delta_0}),\quad |\eta |>1
\end{equation}
for some $\delta_0>0$ independent of $k$.  Using (\ref{h1kdel}) and (\ref{11jun25e8}) in the computation of the second integral of $\rho_i^k$ we have
$$
\rho_i^k=2\pi \bar \sigma_i^k+2\pi \epsilon_k^2 \log \epsilon_k^{-1}
b_i^ke^{c_i^k}+O(\epsilon_k^2).
$$
where
\begin{eqnarray*}
b_i^k=\frac 14\bigg (\frac{\Delta h_i(p_k)}{h_i(p_k)}-2K(p_k)+8\pi+16\pi \frac{\nabla h_i(p_k)}{h_i(p_k)}\nabla_1\gamma(p_k,p_k)
+64\pi^2|\nabla_1 \gamma(p_k,p_k)|^2\bigg )\\
=\frac 14\bigg (\Delta\, \log h_i(p_k)-2K(p_k)+8\pi
+|\nabla \log h_i(p_k)+8\pi \nabla_1\gamma(p_k,p_k)|^2\bigg )
\end{eqnarray*}
Consequently,
\begin{eqnarray}\label{11jun8e1}
&&4\sum_i\frac{\rho_i^k}{2\pi}-\sum_{i,j}a_{ij}\frac{\rho_i^k}{2\pi}\frac{\rho_j^k}{2\pi}\\
&=&\epsilon_k^2\log \epsilon_k^{-1}(\sum_i (4-2m_i^k)b_i^k e^{c_i^k})+O(\epsilon_k^2). \nonumber\\
&=&-4\epsilon_k^2\log \epsilon_k^{-1}\sum_i b_i^k e^{c_i^k} +O(\epsilon_k^2). \nonumber
\end{eqnarray}
Theorem \ref{rhotome4} is established. $\Box$

\begin{rem}
Even though there is a cut-off function in the definition of $\chi$, changing the domain of this cut-off function does not change
$\nabla_1\gamma(p_k,p_k)$.
\end{rem}

\noindent{\bf Proof of Theorem \ref{locations}:}

By going through the proof of Lemma \ref{graest} using Theorem \ref{rhoto1mle4} for $\rho\neq Q$ or Theorem \ref{rhotome4} for $\rho=Q$ instead of Proposition \ref{cruexp} one sees easily that
$$
\phi_i^k(x)=2\pi m_{i}^k (\gamma(p_k+x, p_k)-\gamma(p_k,p_k))
-c_{i}^k-f_{i}^k(x)+O(\epsilon_k^{m^k-2})
$$
if $\rho\neq Q$. On the other hand
$$
\phi_{i}^k(x)=2\pi m_{i}^k (\gamma(p_k+x, p_k)-\gamma(p_k,p_k))
-c_{i}^k-f_{i}^k(x)+O(\epsilon_k^{2}\log \epsilon_k^{-1})
$$
if $\rho=Q$. Correspondingly
\begin{equation}\label{11july13e5}
\nabla \phi_{i}^k(0)=2\pi m_{i}^k \nabla_1\gamma(p_k, p_k)+O(\epsilon_k^{m^k-2})
\end{equation}
if $\rho\neq Q$ and
\begin{equation}\label{11july13e6}
\nabla \phi_{i}^k(0)=2\pi m_{i}^k \nabla_1\gamma (p_k, p_k)+O(\epsilon_k^{2}\log \epsilon_k^{-1})
\end{equation}
if $\rho=Q$. Theorem \ref{location} yields
$$|\sum_i(\partial_l(\log h_{i}^k)(0)+\partial_l \phi_{i}^k(0))\sigma_{i}^k|\le C\epsilon_k^{m^k-2} \quad \mbox{ if } \rho\neq Q $$
and
$$|\sum_i (\partial_l (\log h_{i}^k)(0)+\partial_l\phi_{i}^k(0))\sigma_{i}^k|\le C\epsilon_k^{2}\log \epsilon_k^{-1} \quad \mbox{ if } \rho= Q. $$
Also the proofs of Theorem \ref{rhoto1mle4} and Theorem \ref{rhotome4} give
$$\rho_i^k-\rho_i=O(\epsilon_k^{m^k-2})=O(\epsilon_k^{m-2}), \quad m<4 $$
and
$$\rho_i^k-\rho_i=O(\epsilon_k^{m^k-2})\log \epsilon_k^{-1}=O(\epsilon_k^2\log \epsilon_k^{-1}), \quad m=4.  $$
So
$$\bar \sigma_i^k-\frac{\rho_i}{2\pi}=\left\{\begin{array}{ll} O(\epsilon_k^{m-2}) &\quad m<4,\\
O(\epsilon_k^2\log \epsilon_k^{-1})&\quad m=4.
\end{array}
\right.
$$

By (\ref{11july13e5}) and (\ref{11july13e6}) we obtain (\ref{11july13e7}) and (\ref{11july13e8}) in Theorem \ref{locations}. Theorem \ref{locations} is established $\Box$

\end{document}